\documentclass[10pt]{article}
\usepackage[utf8]{inputenc}
\usepackage{latexsym,amssymb,amsmath,url,authblk}

\def\N{\mathbb{N}}

\def\Z{\mathbb{Z}}
\def\R{\mathbb{R}}
\def\C{\mathbb{C}}

\def\proof{\par\noindent{\em Proof. }}
\def\eproof{\hfill{$\Box$}\bigskip}

\def\ds{\dots}
\def\sus{\subset}
\def\al{\alpha}
\def\be{\beta}

\def\cc{\colon}
\def\ep{\varepsilon}

\newtheorem{thm}{Theorem}[section]
\newtheorem{prop}[thm]{Proposition}
\newtheorem{cor}[thm]{Corollary}

\newtheorem{lem}[thm]{Lemma}

\newtheorem{defi}[thm]{Definition}

\title{Extending P\'olya's random walker beyond probability~I. Complex weights}
\author{M. Klazar\footnote{{\tt klazar@kam.mff.cuni.cz}}\ \ (UK, Praha) and R. Horsk\'y (V\v SE, Praha)}
\date{\today}

\begin{document}

\maketitle

\begin{abstract}
Working in combinatorial model
$\mathrm {W_{co}}(d)$, $d=1,2,\ds$, of P\'olya's
random walker in $\Z^d$, we prove two theorems on 
recurrence to a~vertex. We obtain an effective version of
the first theorem if $d=2$. Using a~semi-formal approach to 
generating functions, we extend both theorems beyond
probability to a~more general model $\mathrm{W_{\C}}$ with
complex weights. We relate models
$\mathrm {W_{co}}(d)$ to standard models 
$\mathrm{W_{Ma}}(d)$ based on Markov chains. The follow-up
article will treat non-Archimedean models $\mathrm{W_{fo}}(k)$ in which weights are formal power series in $\C[[x_1,x_2,\ds,x_k]]$.
\end{abstract}

\tableofcontents

\section[Introduction. P\'olya's walkers]{Introduction. P\'olya's walkers}\label{sec_intro}

In 1921 G.~P\'olya \cite{poly} investigated particles randomly 
moving along edges of the nearest neighbor graph $\Z^d$. We 
call this problem {\em P\'olya's random walker}. It is modeled probabilistically
by Markov chains. In our article we propose and investigate more general 
combinatorial models.

In this Section~\ref{sec_intro} we state in the language of combinatorial 
models $\mathrm{W_{co}}(d)$, $d\in\N$, Theorems~\ref{thm_polya1} 
and \ref{thm_polya2}.
They concern P\'olya's random walker and determine probabilities
of recurrence to a~vertex $\overline{v}$ of walks in $\Z^d$
starting at $\overline{0}$. 
In the former theorem we have $\overline{v}=\overline{0}$, and in 
the latter $\overline{v}\ne\overline{0}$. 
Models $\mathrm{W_{co}}(d)$ are elements of the
more general model $\mathrm{W_{\C}}$ in Definition~\ref{def_WC}. This 
model goes beyond probability\,---\,it uses edge weights in $\C$. In 
conclusion we recall standard models $\mathrm{W_{Ma}}(d)$ based on Markov chains. 

In Section~\ref{sec_relating} we relate $\mathrm{W_{Ma}}(d)$ to 
$\mathrm{W_{co}}(d)$. 
Section~\ref{sec_polya123} contains proofs of Theorems~\ref{thm_polya1} 
and \ref{thm_polya2} in stronger form. In Section~\ref{sec_effeValu} 
we work out an effective version of
Theorem~\ref{thm_polya1} for $d=2$. Using our semi-formal approach to 
generating functions, which extends the symbolic method in combinatorics 
\cite{albe_al,comt,flaj_sedg,goul_jack,symb} from finite to 
countable sets, in Section~\ref{sec_weights}
we extend (strengthened) Theorems~\ref{thm_polya1} and 
\ref{thm_polya2} to $\mathrm{W_{\C}}$. 
Section~\ref{sec_conclR} outlines the follow-up article which will deal with non-Archimedean weights. 
At the very end of this section, we preview our results in more detail.

G.~P\'olya extended his pioneering article \cite{poly} 
in \cite{poly_38}; see \cite{alex} for his 
scientific biography. P\'olya's random walker in $\Z^d$ is discussed, for 
example, in 
\cite{beck,bend_rich,bill,fell,fost_good,grim_wels,koch,lang,levi_pere,
nova,reny,reve,wins,woes}. 
 
We review notation. $\N=
\{1,2,\ds\}$, $\N_0=\N\cup
\{0\}$ and $\Z$ are the integers. By $\R$ and $\C$ we denote
the fields of real and complex numbers, and $\R_{\ge0}$ are 
nonnegative real numbers. 
For a~set $A$ we
denote by $\binom{A}{2}$ the set of two-element subsets of $A$. Instead 
of $\binom{\N}{2}$ we write $\N_2$. If $A$ is a~finite  set, $|A|$ 
($\in\N_0$) denotes the number of its elements. 
$\mathcal{P}(A)
=\{X\cc\;X\sus A\}$ is the power set of $A$. For $n\in\N$ let 
$[n]=\{1,2,\ds,n\}$, and let $[0]=\emptyset$. If $f\cc X\to 
Y$ is a~map and $Z$ is any set, then 
$f[Z]=\{f(x)\cc\;x\in X\cap Z\}$ ($\sus Y$) and 
$f^{-1}[Z]=\{x\in X\cc\;f(x)\in Z\}$ ($\sus X$).

Let $d\in\N$. The nearest neighbor 
graph on $\Z^d$ is the countable graph
$${\textstyle
G_d=(\Z^d,\,E_d),\ E_d\sus\binom{\Z^d}{2}\,,
}
$$
with the vertices
$\Z^d=\Z\times\Z\times\ds\times\Z$ ($d$ factors) and edges 
$$
{\textstyle
\{\overline{u},\,\overline{v}\}=
\{\langle u_1,\,u_2,\,\ds,\,u_d\rangle,\,\langle v_1,\,v_2,\,\ds,\,v_d\rangle\}\in E_d\iff
\sum_{i=1}^d|u_i-v_i|=1\,.
}
$$
$G_d$ has no multiple edges and no loops, and is
$2d$-regular (every vertex lies in $2d$ edges). It is vertex-transitive. 

We consider walks in $G_d$ starting at the vertex $\overline{0}=\langle 0,0\ds,0\rangle$. A~{\em walk} in $G_d$ is any tuple 
$$
w=\langle \overline{u}_0,\,\overline{u}_1,\,\ds,\,\overline{u}_n\rangle,\ n\in\N_0\,,
$$
of vertices $\overline{u}_i\in\Z^d$ such that $\{\overline{u}_{i-1},\overline{u}_i\}\in 
E_d$ for every $i\in[n]$. The set $\{\overline{u}_1,\overline{u}_2,\ds, \overline{u}_{n-1}\}$ is the {\em inner} of $w$. The {\em length $n$} of $w$ is denoted by 
$|w|$. We say that {\em $w$ starts}
at $\overline{u}_0$, or that it is a~{\em $\overline{u}_0$-walk}, and that it {\em ends at $\overline{u}_n$}. We call $w$ {\em
recurrent} if it starts at $\overline{0}$, 
$|w|=n>0$ and $\overline{u}_i=\overline{0}$ for
some $i\in[n]$. Thus $w$ revisits its start $\overline{0}$. Due to the vertex-transitivity of $G_d$, the starting vertex
$\overline{0}$ may
be replaced with any other vertex. For $d\in\N$ and $n\in\N_0$ we let $W_d(n)$ be the
set of $\overline{0}$-walks of length~$n$ in $G_d$. Thus $|W_d(n)|=
(2d)^n$. We denote by $W_{d,\mathrm{rec}}(n)$ ($\sus
W_d(n)$) the set of recurrent walks of length $n$. The next theorem is 
a~combinatorial version of the classical result on random walks 
in $G_d$. It is not P\'olya's theorem which comes in
Theorem~\ref{thm_polya2}. 

\begin{thm}[$\overline{0}$-recurrence]\label{thm_polya1}
For every $d\in\N$, 
$$
\lim_{n\to\infty}\frac{|W_{d,\mathrm{rec}}(\,n)|}{|W_d(n)|}=\lim_{n\to\infty}\frac{|W_{d,\mathrm{rec}}(\,n)|}{(2d)^n}\left\{
\begin{array}{lll}
=1  & \ds  & \text{if $d\le2$ and}  \\
<1  &  \ds & \text{if $d\ge3$}\;.
\end{array}
\right.
$$
\end{thm}
So for $d\le2$ and $n\to\infty$ the proportion of non-recurrent walks goes to $0$. For $d\ge3$ and every
$n\in\N_0$ it is greater than $c>0$. The strengthening in Theorem~\ref{thm_explVal1d3} 
gives a~formula for the limit if $d\ge3$. 

For $d\in\N$, $\overline{v}\in\Z^d$ 
and $n\in\N_0$ let $W_d(\overline{v},n)$ ($\sus 
W_d(n)$) be the set of $\overline{0}$-walks of length 
$n$ in $G_d$ such that $\overline{u}_i=\overline{v}$ for 
some $i\in[n]$. So $W_{d,\mathrm{rec}}(n)=W_d(\overline{0},n)$. We call 
walks in $W_d(\overline{v},n)$ {\em $\overline{v}$-recurrent}. The 
next theorem is a~combinatorial 
version of P\'olya's classical result in \cite{poly}.

\begin{thm}[$\overline{v}$-recurrence]\label{thm_polya2}
For every $d\in\N$ and every $\overline{v}\in\Z^d$,
$$
\lim_{n\to\infty}\frac{|W_d(\overline{v},\,n)|}{|W_d(n)|}=\lim_{n\to\infty}\frac{|W_d(\overline{v},
\,n)|}{(2d)^n}\left\{
\begin{array}{lll}
=1  & \ds  & \text{if $d\le2$ and}  \\
<1  &  \ds & \text{if $d\ge3$}\;.
\end{array}
\right.
$$
\end{thm}
By \cite[Theorem~1.5.4]{lawl} and our Corollary~\ref{cor_RecuInInfi}, 
for every $d\ge3$ the limit is asymptotic to
$$
\frac{\frac{d}{2}\Gamma(\frac{d}{2}-1)\pi^{-d/2}}{\|\overline{v}\|^{d-2}}\ \ \ \ (\|\overline{v}\|\to+\infty)\,,
$$
where $\Gamma(x)=\int_0^{+\infty}t^{x-1}
\mathrm{e}^{-t}\,\mathrm{d}t$ is the Gamma function 
and 
$${\textstyle
\|\overline{v}\|=\sqrt{v_1^2+v_2^2+\ds+
v_d^2}
}
$$ 
is the Euclidean norm. The strengthening in 
Theorem~\ref{thm_explVal2d3} gives  a~formula for the limit if $d\ge3$.

G.~P\'olya announces his result, Theorem~\ref{thm_polya2}, in 
\cite[pp.~149--150]{poly} as follows. 
\begin{quote}
    For concreteness we imagine that the wandering point starts its journey at the origin of 
    coordinates at the time $t=0$ and that it moves with speed $1$. 
    
    $[\ds]$
    
    Now I treat the wandering in the $d$-dimensional net of lines; a~node with the 
    coordinates $a_1,a_2,\ds,a_d$ is given; in question is the probability that 
    in the time span $0<t\le n$ the wandering point coincides at least once with 
    the given node $a_1,a_2,\ds,a_d$. The probability obviously grows with $n$. 
    A~question arises: does it tend to certainty when $n$ grows without limit?

    {\em Yes} when $d=1$ or $d=2$, {\em no} when $d\ge 3$. In the following, I will justify 
    this answer. [our translation, original emphasize]
\end{quote}
P\'olya's results are interpreted in terms of (modern) Markov 
chains, but this is an anachronism. Although A.~Markov 
first published on them in 
1906, rigorous probability theory was born in \cite{kolm} one decade 
after \cite{poly}.

In order to introduce the model $\mathrm{W_{\C}}$
we review complex series. Let $X$ be an at
most countable 
set and $g\cc X\to\C$ be a~map. We call $g$ a~{\em series (on 
$X$)} and denote it by $\sum_{x\in X}g(x)$. We say that 
{\em $\sum_{x\in X}g(x)$ absolutely converges} if two 
equivalent conditions hold.

\begin{enumerate}
\item There is a~constant $c>0$ such
that for every finite set $Y\sus X$ we have $\sum_{x\in Y}|g(x)|\le c$.
\item If $X$ is countable then for every bijection $f\cc\N\to X$ the limit 
$${\textstyle
s=\lim_{n\to\infty}\sum_{i=1}^n g(f(i))\ \  (\in\C) 
}
$$
exists and is independent of $f$.
\end{enumerate}
We call $s$ the {\em sum} of the series and
denote it also by $\sum_{x\in X}g(x)$. For any finite set $X=
\{x_1,x_2,\ds,x_n\}$ any series on $X$ absolutely converges and the sum is just 
$$
g(x_1)+g(x_2)+\ds+g(x_n)\,. 
$$
If $S=\sum_{x\in X}g(x)$ is a~series and $Y\sus X$, we call
$\sum_{x\in Y}g(x)$ a~{\em subseries of $S$}.
It is easy to see that any subseries of an absolutely convergent series 
absolutely converges. If $X$ is countable, $g\cc 
X\to\R$ and for every bijection $f\cc\N\to X$ we have 
$${\textstyle
\lim_{n\to\infty}\sum_{i=1}^n 
g(f(i))=+\infty\,, 
}
$$
we say that the series $\sum_{x\in X}g(x)$ {\em has 
sum $+\infty$}. If $U(x)=\sum_{n=0}^{\infty}u_nx^n$ is 
a~formal power series in $\C[[x]]$, 
we denote the function 
given by sums of $U(x)$ for $x\in\C$
by $F_U(x)$.

Let $n\in\N_0$ and 
$${\textstyle
K_{\N}=(\N,\,\binom{\N}{2})=(\N,\,\N_2)
}
$$ 
be the complete countable 
graph. A~{\em walk} with {\em length} $n$ in $K_{\N}$ is any tuple
$$
w=\langle u_0,\,u_1,\,\ds,\,u_n\rangle,\ n\in\N_0\,,
$$
of vertices $u_i\in\N$ such that $u_{i-1}\ne u_i$ for every $i\in[n]$.
We denote the length $n$ of $w$ by $|w|$. We say that {\em $w$ starts} at $u_0$, 
or that it is a~{\em $u_0$-walk}, and that it {\em ends} at $u_n$. The 
set $\{u_1,u_2,\ds,u_{n-1}\}$ is the {\em inner} of $w$. We 
consider walks starting at $1$. For $n\in\N_0$, let
$$
W(n)=\{w\cc\;\text{$w$ is a~$1$-walk in $K_{\N}$ with $|w|=n$}\}\;.
$$

\begin{defi}[$G(h)$]\label{def_Vh}
Let $X$ be a~set with a~special element $0_X\in X$ and let 
$h\cc\N_2\to X$. We define the graph 
$$
G(h)=(V(h),\,E(h))
$$ 
as the connected component containing $1$ of the subgraph of $K_{\N}$ formed by the edges $e$ with the weight $h(e)\ne 0_X$.    
\end{defi}

If $h\cc\N_2\to X$, where 
$$
X_{\mathrm{ri}}=\langle X,\,0_X,\,1_X,\,+,\,\cdot\rangle
$$ 
is a~commutative ring, we extend $h$ 
to any walk $w=\langle u_0,u_1,\ds,u_n\rangle$ in $K_{\N}$ as follows. For $n=0$ we set $h(w)=1_X$. For $n>0$ we set
$${\textstyle
h(w)=\prod_{i=1}^nh(\{u_{i-1},\,u_i\})\,.
}
$$

If all vertex degrees in the graph $G(h)$ are finite, then 
for every $n\in\N_0$ the series $\sum_{w\in W(n)}h(w)$ has only 
finitely many nonzero summands and has a~sum. More generally, we call 
a~weight $h\cc\N_2\to\C$ {\em light} if for every $n\in\N_0$ the series 
$\sum_{w\in W(n)}h(w)$ absolutely converges. The following definition 
is fundamental in our approach. 

\begin{defi}[$\mathrm{W_{\C}}$]\label{def_WC}
The model $\mathrm{W_{\mathrm{\C}}}$ 
is the set of light weights
$h\cc\N_2\to\C$. We call it P\'olya's complex walker.
\end{defi}
In Section~\ref{sec_weights} we extend Theorems~\ref{thm_polya1} and 
\ref{thm_polya2}, in fact Theorems~\ref{thm_explVal1d3} and 
\ref{thm_explVal2d3}, to the model
$\mathrm{W_{\mathrm{\C}}}$ .

Models $\mathrm{W_{\mathrm{co}}}(d)$ are certain elements in
$\mathrm{W_{\C}}$. The proof of the next proposition is clear and is 
omitted. 

\begin{prop}\label{prop_Wco}
Let $d,n,u\in\N$. Suppose that $h\cc\N_2\to\C$ is such that $G(h)$ 
is isomorphic to $G_d$, that always $h(e)\in\{0,\frac{1}{2d}\}$, that 
$1$ corresponds in the isomorphism to
$\overline{0}$, and $u$
to $\overline{v}$. Then we have for walks $w=\langle u_0,u_1,\ds,u_n\rangle$ in $K_{\N}$ that  
$${\textstyle
(2d)^{-n}\cdot|W_d(\overline{v},\,n)|=\sum_{\substack{w\in W(n)\\
\exists i\in[n]:\,u_i=u}}h(w)\,.
}
$$ 
\end{prop}
The proportion or, if you want, probability of
$\overline{v}$-recurrent walks of 
length $n$ in $G_d$ equals to the
total $h$-weight of $u$-recurrent
walks of length $n$ in $K_{\N}$. 

We review the Markov chain model
$\mathrm{W_{Ma}}(d)$ of P\'olya's random walker. Recall 
that a~probability space is a~triple 
$$
P=\langle\Omega,\,\mathcal{A},\,\mathrm{Pr}\rangle
$$ 
of a~nonempty set $\Omega$,
a~$\sigma$-algebra $\mathcal{A}$ ($\sus\mathcal{P}(\Omega)$) on 
the set $\Omega$ and a~$\sigma$-additive probability map
$\mathrm{Pr}\cc\mathcal{A}\to[0,1]$
such that $\mathrm{Pr}(\Omega)=1$. In more details,  $\mathcal{A}\ne\emptyset$, 
$\mathcal{A}$ is closed to complements to $\Omega$ and to countable unions, and
$$
{\textstyle
\mathrm{Pr}\big(\bigcup_{n=1}^{\infty}A_n\big)=\sum_{n=1}^{\infty}\mathrm{Pr}(A_n)
}
$$
for any mutually disjoint sets $A_n\in\mathcal{A}$, $n\in\N$. Let 
$A,B\in\mathcal{A}$ with
$\mathrm{Pr}(B)>0$. The conditional probability of $A$ given $B$ is the ratio
$${\textstyle
\mathrm{Pr}(A\,|\,B)=\frac{\mathrm{Pr}(A\cap B)}{\mathrm{Pr}(B)}\,.
}
$$
For $\mathrm{Pr}(B)=0$ 
this quantity is not defined. Let $S$ be an at most
countable set. An $S$-valued 
$P$-random variable is any map
$X\cc\Omega\to S$ such that $X^{-1}[\{s\}]\in\mathcal{A}$ for
every $s\in S$. We denote these elements of $\mathcal{A}$ by
``$X=s$''. Their intersections and unions are indicated by $\wedge$ and $\vee$, 
respectively. 

\begin{defi}[$\mathrm{W_{Ma}}(d)$]\label{def_WMa}
Let $P=\langle\Omega,
\mathcal{A},\mathrm{Pr}\rangle$ be a~probability space and 
$$
\mathrm{W_{Ma}}(d)=(X_0,\, X_1,\,X_2,\,\,\ds)
$$
be an infinite sequence
of $\Z^d$-valued $P$-random variables satisfying the following three conditions.  
\begin{enumerate} 
\item  For every $n\in\N_0$ and every $\overline{u}_0$, 
$\overline{u}_1$, $\ds$, $\overline{u}_{n+1}$ in $\Z^d$ with 
$\mathrm{Pr}(X_0=\overline{u}_0
\wedge X_1=\overline{u}_1
\wedge\ds\wedge X_n=\overline{u}_n)>0$ we have
\begin{eqnarray*}
&&\mathrm{Pr}(X_{n+1}=\overline{u}_{n+1}\,|\,X_0=\overline{u}_0\wedge X_1=\overline{u}_1
\wedge\ds\wedge X_n=\overline{u}_n)\\
&&=\mathrm{Pr}(X_{n+1}=\overline{u}_{n+1}\,|\,X_n=\overline{u}_n)\,.
\end{eqnarray*}
\item $\mathrm{Pr}(X_0=\overline{0})=1$.
\item If $n\in\N$,  $\overline{u},\overline{v}\in\Z^d$ and $\mathrm{Pr}(X_{n-1}=\overline{u})>0$, then
$$
\mathrm{Pr}(X_n=\overline{v}\,|\, X_{n-1}=\overline{u})
=\left\{
\begin{array}{lll}
 {\textstyle\frac{1}{2d}}  & \ds & \{\overline{u},\,\overline{v}\}\in E_d\text{ and} \\
  0 & \ds &\{\overline{u},\,\overline{v}\}\not\in E_d\,.
\end{array}
\right.
$$
\end{enumerate}
We call
$\mathrm{W_{Ma}}(d)$ P\'olya's Markov walker
\end{defi}

\centerline{{\bf A~preview of this article}}

\medskip\noindent
In Section~\ref{sec_relating} we obtain six results on relations between models 
$\mathrm{W_{Ma}}(d)$ and
$\mathrm{W_{co}}(d)$. For example,  Corollary~\ref{cor_RecuInInfi} says that for every $d\in\N$ and $\overline{v}\in\Z^d$,
$$
{\textstyle
\mathrm{Pr}\big(\bigvee_{n=1}^{\infty}X_n=\overline{v}\big)=\lim_{n\to\infty}(2d)^{-n}\cdot |W_d(\overline{v},\,n)|\,.
}
$$
In Section~\ref{sec_polya123} we strengthen Theorems~\ref{thm_polya1} and 
\ref{thm_polya2} to
Theorems~\ref{thm_explVal1d3} and 
\ref{thm_explVal2d3}, respectively. The strengthened theorems give
explicit formulas for the limits if $d\ge3$. They are  
$${\textstyle
1-\frac{1}{B(1)}\,\text{ and }\,\sqrt{1-\frac{B_0(1)}{B(1)}}\,, 
}
$$
respectively, where $B(x)$ and $B_0(x)$ are generating functions of
closed walks. In Section~\ref{sec_effeValu} we obtain an effective form of 
Theorem~\ref{thm_polya1} for $d=2$. For example,  
Corollary~\ref{cor_poslCOr} says that for every $N\ge N_0$, with $N_0$ effectively computable, we have
$$
(0.9\log N)^{-1}\le
1-4^{-N}|W_{2,\mathrm{rec}}(N)|
\le(0.1\log N)^{-1}\,.
$$ 
Section~\ref{sec_weights} on  $\mathrm{W_{\C}}$ contains our main results. We 
generalize Theorem~\ref{thm_explVal1d3} to
Theorems~\ref{thm_polya1gen1} and \ref{thm_polya1gen1apul}, and 
Theorem~\ref{thm_explVal2d3} to Theorems~\ref{thm_polya1gen2n} and \ref{thm_polya1gen2vC}. The first theorem
in each pair allows any 
light (and $v$-transitive) weight 
$h\cc\N_2\to\C$
and the computed main quantities
are 
$${\textstyle
\lim_{n\to\infty}\sum_{j=0}^n a_j^h\,\text{ and }\,\lim_{n\to\infty}
\sum_{j=0}^n(a_j^h)'\,,
}
$$
where
$a_j^h$ and $(a_j^h)'$ are, respectively, the total $h$-weights 
of recurrent and $v$-recurrent walks of length 
$j$ in $K_{\N}$. The second theorem in each pair allows only
convex light (and $v$-transitive) weight $h$, but the computed main quantities are the more precise
$${\textstyle
\lim_{n\to\infty}a_n^h\,\text{ and }\,\lim_{n\to\infty}(a_n^h)'\,,
}
$$
respectively.
Convexity of $h$ means roughly that
always $\sum_{v\in\N\setminus\{u\}}h(\{u,v\})=1$.

\section[Relations between models $\mathrm{W_{Ma}}(d)$ and 
$\mathrm{W_{co}}(d)$]{Relations between models $\mathrm{W_{Ma}}(d)$ and $\mathrm{W_{co}}(d)$}\label{sec_relating} 

We establish several relations between models 
$\mathrm{W_{Ma}}(d)$ 
and $\mathrm{W_{co}}(d)$. The next main theorem (which we could not find in the literature) shows 
that P\'olya's Markov walker
correctly recognizes walks among tuples of vertices in the graph $G_d$.

\begin{thm}\label{thm_WcoandWMa}
Let $\mathrm{Pr}$ 
and $X_n$ be as in Definition~\ref{def_WMa}.
If $d\in\N$, $n\in\N_0$ and
$t=\langle\overline{u}_0,\overline{u}_1,
\ds,\overline{u}_n\rangle\in(\Z^d)^{n+1}$, then
$$
\mathrm{Pr}(X_0=\overline{u}_0\wedge X_1=\overline{u}_1\wedge\ds\wedge X_n=\overline{u}_n)=
\left\{
\begin{array}{lll}
{\textstyle\big(\frac{1}{2d}\big)^n}  & \ds  &t\in W_d(n)\text{ and}\\
0  & \ds  &t\not\in W_d(n)\,.
\end{array}
\right.
$$ 
\end{thm}
\proof
Let $n\in\N_0$ and $w=\langle\overline{u}_0,\overline{u}_1,\ds,\overline{u}_n\rangle$ be any $\overline{0}$-walk in $G_d$. 
We prove by induction on 
$i=0,1,\ds,n$ that 
$$
\mathrm{Pr}(X_0=\overline{u}_0\wedge X_1=\overline{u}_1\wedge\ds\wedge X_i=\overline{u}_i)=(2d)^{-i}\,.
$$
For $i=0$ this holds by part~2 of Definition~\ref{def_WMa}. Let $0<i<n$. Then the probability 
$\mathrm{Pr}(X_0=\overline{u}_0\wedge X_1=\overline{u}_1\wedge\ds\wedge X_{i+1}=\overline{u}_{i+1})$ indeed equals
\begin{eqnarray*}
&&\mathrm{Pr}(X_{i+1}=\overline{u}_{i+1}\,|\,X_0=\overline{u}_0\wedge\ds\wedge 
X_i=\overline{u}_i)\cdot\mathrm{Pr}(X_0=\overline{u}_0\wedge\ds\wedge 
X_i=\overline{u}_i)\\
&&=\mathrm{Pr}(X_{i+1}=\overline{u}_{i+1}\,|\,X_i=\overline{u}_i)\cdot(2d)^{-i}\\
&&=(2d)^{-1}\cdot(2d)^{-i}=(2d)^{-i-1}\;.
\end{eqnarray*}
In the first line we use the formula for the conditional
probability and the inductive assumption that
$$
\mathrm{Pr}(X_0=\overline{u}_0\wedge X_1=\overline{u}_1\wedge\ds\wedge X_i=\overline{u}_i)=(2d)^{-i}>0\,.
$$
We use this assumption also in the second line, together with part~1 of Definition~\ref{def_WMa}. 
In the third line we use part~3 of 
Definition~\ref{def_WMa} and the fact that $\{\overline{u}_i,
\overline{u}_{i+1}\}\in E_d$. 

Let $n\in\N_0$ and $\langle\overline{u}_0,\overline{u}_1,\ds,\overline{u}_n\rangle$ 
be in $(\Z^d)^{n+1}\setminus W_d(n)$. Let $j\in[n-1]\cup\{0\}$ be maximum such that 
$$
\langle\overline{u}_0,\,\overline{u}_1,\,\ds,\,\overline{u}_j\rangle\in W_d(j)\;.
$$
If $j$ does not exist, that is if $\overline{u}_0\ne\overline{0}$, we set $j=-1$. By the monotonicity of probability,
$$
\mathrm{Pr}(X_0=\overline{u}_0\wedge\ds\wedge X_n=\overline{u}_n)
\le\mathrm{Pr}(X_0=\overline{u}_0\wedge\ds\wedge X_{j+1}=\overline{u}_{j+1})\;.
$$
If $j=-1$, the last probability is $\mathrm{Pr}(X_0=\overline{u}_0)=0$ 
by part~2 of Definition~\ref{def_WMa}
because $\overline{u}_0\ne\overline{0}$. For $j\ge0$ we set $C=(X_0=\overline{u}_0\wedge X_1=\overline{u}_1\wedge\ds
\wedge X_j=\overline{u}_j)$ ($\in\mathcal{A}$) and get that
\begin{eqnarray*}
\mathrm{Pr}(X_0=\overline{u}_0\wedge\ds\wedge X_{j+1}=\overline{u}_{j+1})&=& 
\mathrm{Pr}(X_{j+1}=\overline{u}_{j+1}\,|\,C)\cdot\mathrm{Pr}(C)\\
&=&\mathrm{Pr}(X_{j+1}=\overline{u}_{j+1}\,|\,X_j=\overline{u}_j)\cdot(2d)^{-j}\\
&=&0\cdot(2d)^{-j}=0
\end{eqnarray*}
---\,thus $\mathrm{Pr}(X_0=\overline{u}_0\wedge X_1=\overline{u}_1\wedge\ds\wedge X_n=\overline{u}_n)=0$.
In the first line we use the formula for the conditional probability and the fact proved in the first case 
that $\mathrm{Pr}(C)=(2d)^{-j}>0$. In the second line 
we use this fact too, and part~1 of 
Definition~\ref{def_WMa}. In the third line we use part~3 of 
Definition~\ref{def_WMa} and the fact that $\{\overline{u}_j,
\overline{u}_{j+1}\}\not\in E_d$ by the maximality of $j$. 
\eproof

\begin{cor}\label{cor_fixEndp}
If $d\in\N$, $n\in\N_0$ and $\overline{v}\in\Z^d$, then
$$
\mathrm{Pr}(X_n=\overline{v})=
(2d)^{-n}\cdot|\{w\in W_d(n)\cc\;\text{$w$ ends at $\overline{v}$}\}|\;.
$$   
\end{cor}
\proof
Let $d$, $n$ and $\overline{v}$ be as stated. If $n=0$, the 
equality holds by part~2 of Definition~\ref{def_WMa}: if 
$\overline{v}=\overline{0}$ (respectively 
$\overline{v}\ne\overline{0}$) then both sides are $1$ (respectively 
$0$). If $n\in\N$, we let 
$$
\widehat{u}=
\langle\overline{u}_0,\,\overline{u}_1,\,\ds,\,\overline{u}_{n-1}\rangle
$$
run in $(\Z^d)^n$ and set $C=(X_0=\overline{u}_0\wedge X_1=\overline{u}_1\wedge\ds\wedge X_{n-1}=\overline{u}_{n-1})$ ($\in\mathcal{A}$). The $\sigma$-additivity of $\mathrm{Pr}(\cdot)$ and Theorem~\ref{thm_WcoandWMa} give 
$${\textstyle
\mathrm{Pr}(X_n=\overline{v})=   \sum_{\widehat{u}}\mathrm{Pr}(C\wedge X_n=\overline{v})
=\sum_{\widehat{u},\;\widehat{u}\,\overline{v}\in W_d(n)}(2d)^{-n}\,,
}
$$
which equals to the stated product.
\eproof

For $\overline{u}=\langle u_1,u_2,\ds,u_d\rangle$ in $\Z^d$ we set $\|\overline{u}\|=
|u_1|+|u_2|+\ds+|u_d|$.

\begin{cor}\label{cor_proNonz}
Let $d\in\N$, $m,n\in\N_0$ and $\overline{u}\in\Z^d$ with $\|\overline{u}\|=n$. Then
$\mathrm{Pr}(X_n=\overline{u})>0$ and
$\mathrm{Pr}(X_m=\overline{u})=0$
for $m<n$.
\end{cor}
\proof
Let $\overline{u}\in\Z^d$ with
$\|\overline{u}\|=n$ ($\in\N_0$). 
Corollary~\ref{cor_fixEndp} gives
$\mathrm{Pr}(X_n=\overline{u})>0$ because there exists a~$\overline{0}$-walk 
$w=\langle\overline{u}_0,\overline{u}_1,\ds,
\overline{u}_n\rangle$ in $G_d$ with $\overline{u}_n=\overline{u}$.
On the other hand, if $m\in\N_0$ and $m<n=\|\overline{u}\|$, then
there is no $\overline{0}$-walk of length $m$ in $G_d$ ending at
$\overline{u}$, and Corollary~\ref{cor_fixEndp} gives 
$\mathrm{Pr}(X_m=\overline{u})=0$.
\eproof

\noindent
We leave it as an exercise for the reader to show that for $m>n=\|\overline{u}\|$ one has
$\mathrm{Pr}(X_m=\overline{u})>0$ iff
$m\equiv n$ modulo $2$.

Recall that $W_d(\overline{v},n)$ is the set of $\overline{0}$-walks 
of length $n$ in the graph $G_d$ such that they visit the vertex $\overline{v}$ at a~step 
$i>0$.

\begin{cor}\label{cor_genRecu}
If $d\in\N$, $n\in\N_0$ and $\overline{v}\in\Z^d$, then
$$
{\textstyle
\mathrm{Pr}\big(\bigvee_{i=1}^n X_i=\overline{v}\big)=(2d)^{-n}\cdot|W_d(\overline{v},\,n)|\,.
}
$$
\end{cor}
\proof
Let $d$, $n$ and $\overline{v}$ be as stated. For $n\in\N_0$ we
let $\widehat{u}=
\langle\overline{u}_0,\overline{u}_1,\ds,\overline{u}_n\rangle$
run in $(\Z^d)^n$. By the $\sigma$-additivity of $\mathrm{Pr}(\cdot)$ 
and Theorem~\ref{thm_WcoandWMa}, the stated probability equals to
$${\textstyle
\sum_{\substack{\widehat{u}\in(\Z^d)^n\\\exists\,i\in[n]:\,\overline{u}_i=\overline{v}}}\mathrm{Pr}(X_0=\overline{u}_0\wedge X_1=\overline{u}_1\wedge\ds\wedge X_n=\overline{u}_n)=
\sum_{\substack{\widehat{u}\in W_d(n)\\\exists\,i\in[n]:\,\overline{u}_i=\overline{v}}}(2d)^{-n}\,,
}
$$
which equals to the stated product.
\eproof

\begin{cor}\label{cor_RecuInInfi}
If $d\in\N$ and $\overline{v}\in\Z^d$, then
$$
{\textstyle
\mathrm{Pr}\big(\bigvee_{n=1}^{\infty}X_n=\overline{v}\big)=\lim_{n\to\infty}(2d)^{-n}\cdot |W_d(\overline{v},\,n)|\,.
}
$$
\end{cor}
\proof
Let $\mathcal{A}$ be as in Definition~\ref{def_WMa}. For every 
sequence of sets $A_n\in\mathcal{A}$, 
$n\in\N$, we have
$$
{\textstyle
\mathrm{Pr}\big(\bigcup_{i=1}^{\infty}A_i\big)=
\lim_{n\to\infty}
\mathrm{Pr}\big(\bigcup_{i=1}^n A_i\big)\,.
}
$$
This follows from the identity 
$${\textstyle
\bigcup_{i=1}^{\infty}A_i=\bigcup_{n
=1}^{\infty}(\bigcup_{i=0}^n A_i\setminus\bigcup_{i=0}^{n-1} 
A_i)
}
$$
(this is a~disjoint union)
where $A_0=\emptyset$, from the equality $\mathrm{Pr}(A\setminus B)=\mathrm{Pr}(A)-\mathrm{Pr}(B)$ 
for sets $B\sus A$ in $\mathcal{A}$, and from the
$\sigma$-additivity of $\mathrm{Pr}(\cdot)$. Since by 
Corollary~\ref{cor_genRecu} we have for every $n\in\N$ that
$${\textstyle
\mathrm{Pr}\big(\bigvee_{i=1}^n X_i=\overline{v}\big)=
(2d)^{-n}\cdot|W_d(\overline{v},\,n)|\,,
}
$$
it follows that $\mathrm{Pr}\big(\bigvee_{n=1}^{\infty}X_n=\overline{v}\big)$ equals
$${\textstyle
\lim_{n\to\infty}\mathrm{Pr}\big(\bigvee_{i=1}^n X_i=\overline{v}\big)=
\lim_{n\to\infty}(2d)^{-n}\cdot|W_d(\overline{v},\,n)|
\,.
}
$$
\eproof

We conclude this section with a~trivial situation which is of some 
interest. For any walk $w=\langle\overline{u}_0,\overline{u}_1,\ds,\overline{u}_n\rangle$
in $G_d$ we denote the last vertex
$\overline{u}_n$ by $\ell(w)$. 

\begin{prop}\label{prop_veryTrivi}
If $d\in\N$, $n\in\N_0$ and $\overline{u},\overline{v}\in\Z^d$ 
with $\overline{u}\ne\overline{v}$, then 
$$
\mathrm{Pr}(X_n=\overline{u}\wedge X_n=\overline{v})=0=(2d)^{-n}\cdot
|\{w\in W_d(n)\cc\;\ell(w)=\overline{u}\wedge\ell(w)=
\overline{v}\}|\,.
$$
\end{prop}
\proof
The latter equality to $0$ is obvious. But so is the former, for $\overline{u}\ne\overline{v}$ we have
$$
X_n^{-1}[\{\overline{u}\}]\cap X_n^{-1}[\{\overline{v}\}]=\emptyset\,.
$$
\eproof

Today the only model of P\'olya's
random walker is 
$\mathrm{W_{Ma}}(d)$. If
$\mathrm{W_{co}}(d)$ is mentioned, it is only as an illustration or motivation. In \cite{klaz} we hope to generalize
the above relations between $\mathrm{W_{Ma}}(d)$ and $\mathrm{W_{co}}(d)$. 
Definition~\ref{def_WMa} is non-vacuous only if the Markov chain
$\mathrm{W_{Ma}}(d)$ and the probability space $P$ supporting it 
exist. They indeed
exist but their rigorous construction requires some effort, see
\cite{bill,woes_MCh}. In contrast, the existence of sets of walks 
$W_d(n)$ and 
$W_d(\overline{v},n)$ in the model $\mathrm{W_{co}}(d)$ is obvious and does not need complicated justification.

\section[Theorems~\ref{thm_polya1} and \ref{thm_polya2} in stronger forms]{Theorems~\ref{thm_polya1} and \ref{thm_polya2} in stronger forms}\label{sec_polya123}

We strengthen
Theorem~\ref{thm_polya1} on $\overline{0}$-recurrence to 
Theorem~\ref{thm_explVal1d3}, and Theorem~\ref{thm_polya2} on 
$\overline{v}$-recurrence to 
Theorem~\ref{thm_explVal2d3}. In the strengthening we obtain 
explicit formulas for the limits if $d\ge3$ (for $d\le2$ the limit is $1$). 

We strengthen Theorem~\ref{thm_polya1}. We
define four sequences $(a_n)$, $(b_n)$, $(c_n)$ and
$(d_n)$ in $\N_0$ counting
certain $\overline{0}$-walks of length $n$ in $G_d$. Let
$n\in\N_0$. We define $a_n$ as the number of recurrent walks, $b_n$ as the number
of walks ending at $\overline{0}$, $c_n$ as the number of walks ending at $\overline{0}$ but with the inner
avoiding $\overline{0}$, and $d_n=(2d)^n$ as the number of all
walks. We have $a_0=c_0=0$ and $b_0=d_0=1$. Let the formal power series
\begin{eqnarray*}
&&{\textstyle
A(x)=\sum_{n=0}^{\infty}\frac{a_n}{d_n}x^n,\  
B(x)=\sum_{n=0}^{\infty}\frac{b_n}{d_n}x^n,\  
C(x)=\sum_{n=0}^{\infty}\frac{c_n}{d_n}x^n}\\ 
&&{\textstyle \text{and }\, D(x)=\sum_{n=0}^{\infty}\frac{d_n}{d_n}x^n=\frac{1}{1-x}
}    
\end{eqnarray*}
in $\C[[x]]$ be the corresponding generating functions. The numbers $d_n$ are multiplicative: for every $m,n\in\N_0$ we have $d_{m+n}=d_m d_n$. 
We suppress in the notation the dependence on $d$. The quantities of
the primary interest are $a_n$ and $A(x)$. 

The above generating functions satisfy two relations. 

\begin{prop}\label{prop_firstRel}
Let $A(x)$, $B(x)$, $C(x)$ and $D(x)$ be as above. The following formal relations hold. 
\begin{enumerate}
\item  
$A(x)=C(x)D(x)=C(x)\frac{1}{1-x}$.
\item $B(x)=\frac{1}{1-C(x)}=\sum_{j=0}^{\infty}C(x)^j$, equivalently
$C(x)=1-\frac{1}{B(x)}$.
\end{enumerate}  
\end{prop}
\proof
1. If $n\in\N_0$ then
$$
{\textstyle
a_n=\sum_{j=0}^n c_jd_{n-j},
\,\text{ so that }\,\frac{a_n}{d_n}=\sum_{j=0}^n \frac{c_j}{d_j}\cdot1\,.
}
$$
The first equality follows from the unique  splitting of any recurrent 
walk $w$ of length $n$ in $G_d$ at the first revisit of 
$\overline{0}$ in step $j$ as $w=w_1w_2$. The walks $w_1$ are counted by $c_j$, and the walks $w_2$ by $d_{n-j}$. Dividing by $d_n$ and using that $d_n=d_jd_{n-j}$ we get
the second equality. In $\C[[x]]$
it reads as $A(x)=C(x)D(x)$.

2. Now $b_0=1$ and for $n\in\N$,
$${\textstyle
b_n=\sum_{j=1}^{\infty}\sum_{\substack{l_1,\ds,l_j\in\N\\l_1+\ds+l_j=n}}c_{l_1}c_{l_2}\ds c_{l_j}\,.
}
$$
This equality follows from the unique decomposition of any 
$\overline{0}$-walk of length $n\in\N$ ending at 
$\overline{0}$ into $j\in\N$ walks with lengths $l_i$ such that each 
starts and ends at $\overline{0}$ 
and has inner avoiding $\overline{0}$. Dividing by $d_n$ and using the multiplicativity of $d_n$ we get 
that for every $n\in\N_0$,
$${\textstyle
\frac{b_n}{d_n}=\sum_{j=0}^{\infty}\sum_{\substack{l_1,\ds,l_j\in\N\\l_1+\ds+l_j=n}}\frac{c_{l_1}}{d_{l_1}}\cdot\frac{c_{l_2}}{d_{l_2}}\cdot\ldots\cdot\frac{c_{l_j}}{d_{l_j}}\;.
}
$$
Multiplying it by $x^n$ and formally summing the results over $n\in\N_0$ we get the relation
$B(x)=\sum_{j=0}^{\infty}C(x)^j=\frac{1}{1-C(x)}$. Solving it for $C(x)$ 
we get the latter equality.
\eproof

We denote by $\R_{\ge0}[[x]]$ the univariate formal power series with non-negative real coefficients.

\begin{defi}[$U(R)$ 1]\label{def_UzR}
Let $R>0$ be a~real number and  $U(x)=\sum_{n=0}^{\infty}u_nx^n$ be in $\R_{\ge0}[[x]]$. We define
$${\textstyle
U(R):=\lim_{n\to\infty}\sum_{j=0}^n u_jR^j\ \ (=:
\sum_{n=0}^{\infty}u_nR^n)\,. 
}
$$
$U(R)$ is a~nonnegative real number or $+\infty$.
\end{defi}
This limit always exists because partial sums $\sum_{j=0}^n u_jR^j$, $n=0,1,\ds$, are non-decreasing. 

\begin{prop}\label{prop_tauber1}
Suppose that $U(x)=\sum_{n=0}^{\infty}u_n x^n$ and $V(x)=\sum_{n=0}^{\infty}v_n x^n$ 
are in $\R_{\ge0}[[x]]$ and that formally $U(x)=V(x)\frac{1}{1-x}$.  
Then
$$
\lim_{n\to\infty}u_n=
\lim_{n\to\infty}(v_0+v_1+\ds+v_n)=V(1)\ \ (\in[0,\,+\infty)\cup\{+\infty\})\,.
$$
\end{prop}
\proof
This follows from $u_n=\sum_{j=0}^nv_j$, $n\in\N_0$.
\eproof

For the sake of completeness we prove the next folklore result which 
we call the {\em weak Abel theorem}. The full Abel theorem \cite{abel} comes in Theorem~\ref{thm_abel}.

\begin{thm}[weak Abel]\label{thm_weakAbel}
Let $R>0$ be a~real number and let $U(x)=\sum_{n=0}^{\infty}u_nx^n$ be 
in $\R_{\ge0}[[x]]$. Suppose that for every $x\in[0,R)$ the series 
$U(x)$ absolutely converges and consider its sum $F_U\cc[0,R)\to[0,+\infty)$. 
Then 
$$
\lim_{x\to R}F_U(x)=U(R)\ \ (\in[0,\,+\infty)\cup\{+\infty\})\,.
$$
\end{thm}
\proof
Let $R$ and $U(x)$ be as stated, and let $N\in\N_0$. Since $R$ and all coefficients $u_n$ are nonnegative, 
\begin{eqnarray*}
&&{\textstyle
\sum_{n=0}^N u_nR^n=\lim_{x\to R^-}\sum_{n=0}^N u_nx^n\le
}\\
&&{\textstyle
\le\lim_{x\to R}F_U(x)=\lim_{x\to R}\sum_{n=0}^{\infty} u_nx^n\le
\sum_{n=0}^{\infty} u_nR^n=U(R)\,.
}
\end{eqnarray*}
For $N \to +\infty$ we get
$U(R)=\lim_{x\to R}F_U(x)$.
\eproof 

We obtain two corollaries of the weak Abel theorem.

\begin{cor}\label{cor_UV1}
Suppose that $U(x)$ and $V(x)$ in $\R_{\ge0}[[x]]$
converge on $[0,1)$, $U(0)=u_0>0$ and that formally
$V(x)=1-\frac{1}{U(x)}$. Then $U(1)>0$ and if $U(1)<+\infty$,
$${\textstyle
V(1)=1-\frac{1}{U(1)}\,. 
}
$$
If $U(1)=+\infty$ then $V(1)=1$.
\end{cor}
\proof
Formally
$V(x)U(x)=U(x)-1$. It follows that $F_U>0$ on $[0,1)$ and that $F_U$ is non-decreasing. By Propositions~\ref{prop_ACser1apul}, \ref{prop_grouSer} 
and \ref{prop_proSer} (which we prove later),
$$
F_V(x)F_U(x)=F_U(x)-1
$$
and $F_V(x)=1-\frac{1}{F_U(x)}$ for every $x\in[0,1)$. 
For $x\to 1$ this equality gives by Theorem~\ref{thm_weakAbel} the two claims of the corollary
\eproof

\begin{cor}\label{cor_UVW1}
Suppose that $U(x)$, $V(x)$ 
and $W(x)$ are in $\R_{\ge0}[[x]]$ and converge on $[0,1)$, 
$W(0)=w_0>0$, $W(1)<+\infty$ and
that formally
$W(x)=V(x)+W(x)U(x)^2$. Then 
$W(1)>0$, 
$\frac{V(1)}{W(1)}\le1$ and
$${\textstyle
U(1)=\sqrt{1-\frac{V(1)}{W(1)}}\,.
}
$$
\end{cor}
\proof
$W(1)>0$ is immediate from $w_0>0$. Arguing as in the previous proof we get the equality $W(1)=V(1)+W(1)U(1)^2$. Since $U(1)\ge0$, we have $U(1)=\sqrt{1-\frac{V(1)}{W(1)}}$.
\eproof

\noindent
Alternatively, we can avoid functions $F_U(x)$, $F_V(x)$, 
and $F_W(x)$ and obtain Corollaries~\ref{cor_UV1} and 
\ref{cor_UVW1} by
manipulating only series with nonnegative summands, convergent
ones and the ones with the sum $+\infty$

We review the asymptotic symbols $\ll$, $\gg$, $O$ and $\sim$. Let $f,g\cc 
X\to\C$ where $X$ is a~nonempty set. We write 
$$
\text{$f\ll g$ on $X$} 
$$
if 
there is a~real constant $c>0$ such that $|f(x)|\le c\cdot|g(x)|$ for every $x\in X$.
The notation $f\gg g$ means $g\ll f$, and $f=O(g)$ is 
synonymous with $f\ll g$. If $X\sus\R$ and $A\in\R\cup\{-
\infty,+\infty\}$ is a~limit point of the definition domain of $\frac{f}{g}$, we write 
$$
\text{$f\sim g$ ($x\to A$)} 
$$
if $\lim_{x\to A}\frac{f(x)}{g(x)}=1$. If $X=\N$, 
$A=+\infty$ and  $f\sim g$ ($x\to+\infty$) then both 
$f\ll g$ and $f\gg g$ hold on the definition domain of $\frac{f}{g}$. 

The next 
Propositions~\ref{prop_bnFordtwo} and \ref{prop_bnFordthree}  deal 
with the generating function $B(x)$ defined above. We show that 
$B(1)=+\infty$ for $d\le2$ and $B(1)<+\infty$ for $d\ge3$. These
are well-known results. We prove them here in detail because we want
to have proofs for them that we like. Let $b_n$, $d_n$ and $B(x)$ be 
as defined above. 

\begin{prop}\label{prop_bnFordtwo}
Let $d\le2$. Then $b_n=0$ iff $n\in\N$ is odd, and
$${\textstyle
\frac{b_n}{d_n}=\frac{b_n}{(2d)^n}\gg n^{-d/2}
}
$$
on even $n\in\N_0$. Hence $B(1)=+\infty$ for $d\le2$.
\end{prop}
\proof
The fact that the $b_n$ are supported on even $n$ is true in 
every dimension $d$ because every walk in $G_d$ that starts and ends 
at $\overline{0}$ does in each of the $d$ directions equally many 
positive and negative steps. For even $n$ obviously $b_n\ge1$. As for 
the estimate 
$\gg$, we prove the stronger relation $\sim$.

Let $d=1$ and $n=2m$, $m\in\N_0$. Then $b_n=b_{2m}=
\binom{2m}{m}$. Hence $\frac{b_n}{d_n}=4^{-m}\binom{2m}{m}$. The Stirling formula 
$m\sim\sqrt{2\pi m}
\big(\frac{m}{\mathrm{e}}\big)^m$ 
($m\to+\infty$), see \cite{robb} 
and Theorem~\ref{thm_effStirl}, gives  
$${\textstyle
\frac{b_n}{d_n}=4^{-m}\cdot(2m)!\cdot m!^{-2}\sim\frac{1}{\sqrt{\pi}}\cdot m^{-1/2}=\sqrt{2/\pi}\cdot n^{-1/2} \ \ (n\to+\infty)\,.
}
$$
Let $d=2$. Then 
$${\textstyle
b_n=b_{2m}=\sum_{j=0}^m\frac{(2m)!}{j!(m-j)!j!(m-j)!}=\binom{2m}{m}\sum_{j=0}^m\binom{m}{j}^2=\binom{2m}{m}^2\,.
}
$$
The first summand is the number of 
walks in $G_2$ that start and end at $\overline{0}$ and do $j$ steps left, $j$ steps right, $m-j$ steps up and $m-j$ 
steps down. The last equality follows 
from the well-known 
binomial identity $\sum_{j=0}^m
\binom{m}{j}^2=\binom{2m}{m}$. The asymptotic relation
$${\textstyle
\frac{b_n}{d_n}=\binom{2m}{m}^2\cdot 4^{-2m}\sim
\frac{1}{\pi m}=\frac{2}{\pi}\cdot n^{-1} \ \  (n\to+\infty)
}
$$
follows again from the Stirling formula. By Definition~\ref{def_UzR} and since $\sum_{n=1}^{\infty}n^{-1}=+\infty$, for $d\le2$ we have  $B(1)=+\infty$.
\eproof

\begin{prop}\label{prop_bnFordthree}
Let $d\ge 3$. Then
$${\textstyle
\frac{b_n}{d_n}\ll n^{-d/2}
}
$$
on $\N$. Hence $B(1)<+\infty$ for $d\ge3$.
\end{prop}
\proof
Let $d\ge3$.
As we know, $b_n=0$ for odd $n$. For $n=2m$ with $m\in\N_0$ the fraction $\frac{b_n}{d_n}=\frac{b_{2m}}{d_{2m}}$ equals 
\begin{eqnarray*}
&{\hspace{3mm}}&{\textstyle(2d)^{-2m}\sum_{\substack{m_i\in\N_0\\m_1+\ds+m_d=m}}\frac{(2m)!}{(m_1!\cdot\ldots\cdot m_d!)^2}
=4^{-m}\binom{2m}{m}\sum_{\cdots}\big[\frac{1}{d^m}\binom{m}{m_1,\,\ds,\,m_d}\big]^2}\\
&\le&{\textstyle 4^{-m}\binom{2m}{m}\cdot\max_{\cdots}d^{-m}\binom{m}{m_1,\,\ds,\,m_d}
=4^{-m}\binom{2m}{m}\cdot d^{-m}\binom{m}{m_1',\,\ds,\,m_d'}}\\
&=&{\textstyle\frac{1}{\sqrt{\pi m}}\cdot
\mathrm{e}^{O(1)}\cdot m^{(1-d)/2}\ll n^{-d/2},\ 
\text{$m_i'\in\{\lfloor\frac{m}{d}\rfloor,\lceil\frac{m}{d}\rceil\}$ and $\sum_{i=1}^d m_i'=m$}\,. 
}
\end{eqnarray*}
The numbers of 
$\lfloor\frac{m}{d}\rfloor$s and $\lceil
\frac{m}{d}\rceil$s in the tuple $m_1'$, $m_2'$, $\ds$, $m_d'$ are 
uniquely determined. In the first line we count the walks as in the case $d\le2$. The numbers in $[\ds]$ sum to $1$ due to the
multinomial expansion of 
$d^m=(1+\ds+1)^m$. The inequality in the second line follows
from the lemma that if $a_1,\ds,a_r\ge0$ and
$a_1+\ds+a_r=1$ then $a_1^2+\ds+a_r^2\le\max_{1\le i\le r}a_i$. The 
multinomial coefficient is maximized by the lemma that if $a > b 
\ge 0$ are integers with 
$a \ge b + 2$ then $a!\cdot b! > (a - 1)!\cdot(b + 1)!$. 
In the third line we use in the following way the asymptotics from the proof of Proposition~\ref{prop_bnFordtwo} 
and the Stirling formula. As $m\to+\infty$, 
$${\textstyle
\binom{m}{m_1',\,\ds,\,m_d'}=\frac{m!}{\prod_{i=1}^d m_i'!}\sim
\frac{\sqrt{2\pi m}\cdot(m/\mathrm{e})^m}{\prod_{i=1}^d\sqrt{2\pi m_i'}\cdot(m_i'/\mathrm{e})^{m_i'}}
=\mathrm{e}^{O(1)}\cdot m^{(1-d)/2}\cdot d^m
}
$$
because $\sum_{i=1}^d m_i'=m$, $\sqrt{m_i'}=(1+O(m^{-1}))\sqrt{m/d}$, $\log(m_i')=\log(\frac{m}{d})+O(m^{-1})$ and
$$
{\textstyle
\big(m_i'\big)^{m_i'}=\mathrm{e}^{m_i'\log(m_i')}=\big(\frac{m}{d}\big)^{m_i'}\cdot\mathrm{e}^{O(m_i'/m)}\,.
}
$$
Here $\mathrm{e}^{O(1)}$ denotes a~function $f\cc\N_0\to(0,+\infty)$ 
such that $1\ll f\ll 
1$ on $\N_0$. Since the series $\sum_{n=1}^{\infty}n^{-3/2}$ 
absolutely converges, for $d\ge3$ we have $B(1)<+\infty$.
\eproof

The first main result of this section is a~strengthening of 
Theorem~\ref{thm_polya1} on $\overline{0}$-recurrence.

\begin{thm}[$\overline{0}$-recurrence${}^+$]\label{thm_explVal1d3}
Let $d\in\N$. In the above notation we have 
$$
\lim_{n\to\infty}\frac{a_n}{d_n}=\lim_{n\to\infty}\frac{|W_{d,\mathrm{rec}}(n)|}{(2d)^n}\left\{
\begin{array}{lll}
=1  & \ds  & \text{if $d\le2$ and}  \\
=1-\frac{1}{B(1)}\in(0,\,1)  &  \ds & \text{if $d\ge3$}\,,
\end{array}
\right.
$$
where $B(x)$ is the generating function of closed walks.
\end{thm}
\proof
By part~1 of Proposition~\ref{prop_firstRel},
formally $A(x)=C(x)\frac{1}{1-x}$. By Proposition~\ref{prop_tauber1}, $\lim_{n\to\infty}
\frac{a_n}{d_n}=C(1)$. By part~2 of Proposition~\ref{prop_firstRel} and 
Corollary~\ref{cor_UV1}, 
$C(1)=1-\frac{1}{B(1)}$. If $d\le2$ then by
Proposition~\ref{prop_bnFordtwo} and 
Corollary~\ref{cor_UV1}, $C(1)=1$.
If $d\ge3$ then by Proposition~\ref{prop_bnFordthree} and 
Corollary~\ref{cor_UV1}, $C(1)=1-
\frac{1}{B(1)}\in(0,1)$.
\eproof

We strengthen  Theorem~\ref{thm_polya2}, the theorem G.~P\'olya proved in \cite{poly}. Let $d\in\N$, 
$\overline{v}\in\Z^d\setminus\{\overline{0}\}$ (we have just treated the case $\overline{v}
=\overline{0}$) and 
$n\in\N_0$. We again define numbers of certain 
$\overline{0}$-walks of length $n$ ($\in\N_0$) in $G_d$. Let $a_n'$ be the number 
of $\overline{v}$-recurrent walks, $b_n$ be as before the number 
of walks ending at 
$\overline{0}$, $b_n'$ be the number of walks ending at $\overline{0}$ and avoiding 
$\overline{v}$, $c_n'$ be the number of walks ending at $\overline{v}$ and with the inner 
avoiding $\overline{v}$, $c_n''$ be the number of walks ending at $\overline{v}$ 
and with the inner avoiding both $\overline{0}$ and $\overline{v}$, and $d_n=(2d)^n$ be as before the number of all walks. We have 
$a_0'=c_0'=c_0''=0$ and $b_0=b_0'=d_0=1$. Besides the previous generating functions $B(x)$ and 
$D(x)=\frac{1}{1-x}$ we consider the 
generating functions 
\begin{eqnarray*}
&&{\textstyle
A_0(x)=\sum_{n=0}^{\infty}\frac{a_n'}{d_n}x^n,\  B_0(x)=\sum_{n=0}^{\infty}
\frac{b_n'}{d_n}x^n,\  
C_0(x)=\sum_{n=0}^{\infty}\frac{c_n'}{d_n}x^n}\\ 
&&{\textstyle\text{and } 
C_1(x)=\sum_{n=0}^{\infty}\frac{c_n''}{d_n}x^n\,.
}    
\end{eqnarray*}
The dependence on 
$d$ and $\overline{v}$ is left implicit in the notation. The 
quantities we want to determine are $a_n'$ and 
$A_0(x)$.

We obtain an analog of Proposition~\ref{prop_firstRel}. 

\begin{prop}\label{prop_firstRel2}
Let $A_0(x)$, $B(x)$, $B_0(x)$, $C_0(x)$, $C_1(x)$, 
and $D(x)$ be the above generating functions. The 
following formal relations hold. 
\begin{enumerate}
    \item $A_0(x)=C_0(x)D(x)=C_0(x)\frac{1}{1-x}$.
    \item $B(x)=B_0(x)+C_0(x)^2B(x)$,
    equivalently $B(x)=\frac{B_0(x)}{1-C_0(x)^2}$.
    \item $C_0(x)=B_0(x) C_1(x)$. 
\end{enumerate}
\end{prop}
\proof
1. The walks $w$ counted by $A_0(x)$ split at the 
first visit of $\overline{v}$ as $w=w_1w_2$, where $w_1$
are walks counted by $C_0(x)$ and $w_2$ are arbitrary
$\overline{v}$-walks. By the vertex-transitivity of $G_d$, 
the walks $w_2$ are counted by $D(x)$.

2. It suffices to prove the first equality. 
The generating function $B(x)$ counts 
$\overline{0}$-walks $w$ in $G_d$ ending at 
$\overline{0}$. Those avoiding $\overline{v}$ are counted 
by $B_0(x)$, and those visiting $\overline{v}$ 
decompose as $w=w_1w_2w_3$, where the walk $w_1$ starts at 
$\overline{0}$, ends at $\overline{v}$ and the inner 
avoids $\overline{v}$, the walk $w_2$ starts and 
ends at $\overline{v}$, and the walk $w_3$ starts at 
$\overline{v}$, ends at $\overline{0}$
and the inner avoids $\overline{v}$. Reversing $w_3$ we see that both $w_1$ 
and $w_3$ are counted by $C_0(x)$. By the vertex-transitivity of $G_d$, the middle
walks $w_2$ are counted by $B(x)$.
We get the first equality. 

3. We consider the walks $w$ counted by $C_0(x)$. They 
start at $\overline{0}$, end at 
$\overline{v}$ and the inner avoids $\overline{v}$. Each $w$ splits at the last visit of $\overline{0}$ in two walks as 
$w=w_1w_2$, 
where $w_1$ and $w_2$ are counted by $B_0(x)$ 
and $C_1(x)$, respectively. The relation follows. 
\eproof

The second main result of this section is a~strengthening of 
Theorem~\ref{thm_polya2} on $\overline{v}$-recurrence.

\begin{thm}[$\overline{v}$-recurrence${}^+$]\label{thm_explVal2d3}
Let $d\in\N$ and
$\overline{v}\in\Z^d$, $\overline{v}\ne\overline{0}$. In 
the above notation we have 
$$
\lim_{n\to\infty}\frac{a_n'}{d_n}=\lim_{n\to\infty}\frac{|W_d(\overline{v},
\,n)|}{(2d)^n}\left\{
\begin{array}{lll}
=1  & \ds  & \text{if $d\le2$ and}  \\
=
\sqrt{1-\frac{B_0(1)}{B(1)}}\in(0,\,1)  &  \ds & \text{if $d\ge3$}\,,
\end{array}
\right.
$$
where $B_0(x)$ and $B(x)$ are the generating functions of closed walks.
\end{thm}
\proof
The generating functions $B(x)$, $B_0(x)$, $C_0(x)$ and $C_1(x)$
have all coefficients in
$[0,1]$. Thus their sums are continuous and non-decreasing functions 
$$
F_B,\,F_{B_0},\,F_{C_0},\,F_{C_1}\cc
[0,\,1)\to[0,\,+\infty)\,.
$$
We show that $F_{C_0}$ is bounded. For $C_0(x)=0$ it certainly 
holds. Suppose that
$C_0(x)$ has some positive coefficient. In fact, it has
infinitely many of them, but we give our argument with
minimal assumptions so that we can use it later in a~more general
situation. If $F_{C_0}<1$, we are done. Else, since 
$F_{C_0}(0)=0$ and $F_{C_0}$ is continuous and increasing, for some 
point $y\in(0,1)$ we have 
$F_{C_0}<1$ on $[0,y)$ and
$\lim_{x\to y^-}F_{C_0}(x)=1$. Since $F_{B_0}(0)=1$, by part~2 of 
Proposition~\ref{prop_firstRel2} we have 
$$
\lim_{x\to y^-}
F_{B}(x)=+\infty\,.
$$
This contradicts the continuity of $F_{B}$ at $y$. Hence $F_{C_0}$ is 
bounded. We finally deduce that 
$F_{B_0}$ is bounded. Else part~3 of 
Proposition~\ref{prop_firstRel2} gives, because $C_1(x)$ has at least one positive coefficient, that
$$
\lim_{x\to 1}F_{C_0}(x)=
\lim_{x\to 1}F_{B_0}(x)\cdot \lim_{x\to 1}F_{C_1}(x)=+\infty\,,
$$
in contradiction with the boundedness of $F_{C_0}$

Let $d\le 2$. For $x\to 1$ we get, using the
boundedness of $F_{B_0}$, Proposition~\ref{prop_bnFordtwo}, part~2 of Proposition~\ref{prop_firstRel2} and Theorem~\ref{thm_weakAbel}, that
$C_0(1)=1$. By part~1 of Proposition~\ref{prop_firstRel2} 
and Proposition~\ref{prop_tauber1}, 
$${\textstyle
\lim_{n\to\infty}\frac{a_n'}{d_n}=C_0(1)=1\,.
}
$$

Let $d\ge3$. Part~2 of Proposition~\ref{prop_firstRel2}, 
Proposition~\ref{prop_bnFordthree} and Corollary~\ref{cor_UVW1} give
$C_0(1)=\sqrt{1-\frac{B_0(1)}{B(1)}}$. By part~1 of Proposition~\ref{prop_firstRel2} 
and Proposition~\ref{prop_tauber1} we have
$$
{\textstyle
\lim_{n\to\infty}\frac{a_n'}{d_n}=
C_0(1)=\sqrt{1-\frac{B_0(1)}{B(1)}}\,.
}
$$
Clearly, $0<B_0(1)<B(1)<+\infty$. Thus the root lies in $(0,1)$.
\eproof

\begin{cor}\label{cor_B0lomenoB}
If $d\le2$ then $\lim_{x\to1}F_{B}(x)=+\infty$ and
$\lim_{x\to1}\frac{F_{B_0}(x)}{F_{B}(x)}=0$.
\end{cor}
\proof
Let $d\le2$. The first limit follows from Theorem~\ref{thm_weakAbel} and Proposition~\ref{prop_bnFordtwo}. 
Part~2 of Proposition~\ref{prop_firstRel2}, and Propositions~\ref{prop_ACser1apul}, \ref{prop_grouSer} and
\ref{prop_proSer} imply that for every $x\in[0,1)$,
$${\textstyle
\frac{F_{B_0}(x)}{F_{B}(x)}= 
}1-F_{C_0}(x)^2\,.
$$
In the previous proof we showed  that for $d\le2$,  
$\lim_{x\to1}F_{C_0}(x)=1$. The stated limit follows.
\eproof

If $d\ge3$ and one sets  $\overline{v}=\overline{0}$, then
Theorem~\ref{thm_explVal2d3} does not turn in 
Theorem~\ref{thm_explVal1d3}: $B_0(1)$ turns in $C(1)$ and, by part~2 of Proposition~\ref{prop_firstRel},
$$
{\textstyle
\sqrt{1-\frac{C(1)}{B(1)}}=
\sqrt{\big(1-\frac{1}{B(1)}\big)^2+\frac{1}{B(1)}}>
1-\frac{1}{B(1)}\,.
}
$$

\section[An effective version of Theorem~\ref{thm_polya1} for $d=2$]{Effective Theorem~\ref{thm_polya1} for $d=2$}\label{sec_effeValu}

Let $f\cc\N\to\R$. The bare limit
$$
\lim_{n\to\infty}f(n)=c
$$ 
is an ineffective result, for no $n\in\N$ it tells us anything
about the distance between $f(n)$ and $c$. Theorems~\ref{thm_polya1},  \ref{thm_polya2}, 
\ref{thm_explVal1d3} and \ref{thm_explVal2d3} on 
$\overline{0}$- and $\overline{v}$-recurrence have this ineffective 
form. We improve upon it a~little and in 
Propositions~\ref{prop_Po_1_2_ub} and \ref{prop_Po_1_2_lb} we 
effectivize Theorem~\ref{thm_polya1} for $d=2$. These effective bounds 
could be further generalized, but for brevity we are content with just 
one case. 

We use an effective Stirling formula.

\begin{thm}[Robbins \cite{robb}]\label{thm_effStirl}
For every $n\in\N$,
$${\textstyle
\sqrt{2\pi n}\cdot\big(\frac{n}{\mathrm{e}}\big)^n\cdot\mathrm{e}^{1/(1+12n)}\le n!\le\sqrt{2\pi n}\cdot\big(\frac{n}{\mathrm{e}}\big)^n\cdot
\mathrm{e}^{1/12n}\,.
}
$$    
\end{thm}
For example,
$$
5039.33\ldots\le 7!=5040\le 5040.04\ds\,.
$$

We obtain an effective version of Corollary~\ref{cor_UV1} if $U(1)=+\infty$.

\begin{prop}\label{prop_effWeakAbel1}
Suppose that the formal power series
$${\textstyle
U(x)=\sum_{n\ge0}u_nx^n\, 
\text{ and }\,V(x)=\sum_{n\ge0}v_nx^n
}
$$ 
in $\R_{\ge0}[[x]]$ absolutely converge for every $x\in[0,1)$. 
Let $u_0>0$, $v_n\le1$ and $u_{2n}\ge\frac{c}{n}$ for every $n$ 
in $\N$ and some constant $c>0$, and let formally $V(x)=1-
\frac{1}{U(x)}$. Then 
$V(1)=\sum_{n=0}^{\infty}v_n=1$ 
and for every $N\in\N$
and $x\in(\frac{1}{2},1)$ we have
$$
1\ge\sum_{n=0}^N v_n\ge 1-\frac{1}{c\log\big(\frac{1}{1-x}\big)-c\log 2}-\frac{x^{N+1}}{1-x}\,.
$$
\end{prop}
\proof
Since, by the assumption, $U(1)=+\infty$, Corollary~\ref{cor_UV1} gives
$V(1)=1$. We get the former
displayed inequality
because $u_n\ge0$. We prove the latter.
Since $u_{2n}\ge\frac{c}{n}$, for every $x\in(0,1)$ we have
$${\textstyle
F_U(x)=\sum_{n=0}^{\infty}u_nx^n\ge 
c\sum_{n=1}^{\infty}\frac{x^{2n}}{n}=c\log\big(\frac{1}{1-x^2}\big)\,.
}
$$
Since $v_n\le 1$ and for 
every $x\in[0,1)$ we have (by the proof of Corollary~\ref{cor_UV1}) 
$F_V(x)=1-\frac{1}{F_U(x)}$, for every $N\in\N$ and 
$x\in(\frac{1}{2}, 1)$ we have
\begin{eqnarray*}
&&{\textstyle\sum_{n=0}^N v_n\ge\sum_{n=0}^N v_nx^n=\sum_{n=0}^{\infty} v_nx^n-\sum_{n>N} v_nx^n}\\
&&{\textstyle\ge 1-\frac{1}{F_U(x)}-\frac{x^{N+1}}{1-x}\ge1+\frac{1}{c\log(1-x)+c\log(1+x)}-\frac{x^{N+1}}{1-x}}\\
&&{\textstyle\ge1-\frac{1}{c\log(\frac{1}{1-x})-c\log 2}-\frac{x^{N+1}}{1-x}\,.}
\end{eqnarray*}
\eproof

\noindent

From this we deduce an effective upper bound on the difference
$$
1-4^{-N}|W_{\mathrm{2,rec}}(N)|
$$ 
in Theorem~\ref{thm_polya1}
for $d=2$. 

\begin{prop}\label{prop_Po_1_2_ub}
For every $N\ge3$  we have
$$
0\le 1-4^{-N}|W_{\mathrm{2,rec}}(N)|\le
(0.2\log N - 0.16)^{-1}+N^{8/9}\exp(-N^{1/9})\;.
$$
\end{prop}
\proof
Let $d=2$, $n\in\N_0$ and let $b_n$ and $c_n$ be as defined at the
beginning of Section~\ref{sec_polya123}. We set 
$v_n=c_n4^{-n}$ and $u_n=b_n4^{-n}$. Then $v_n, u_n \ge0$, $v_n\le 1$, $u_0>0$, the power series
$V(x) = C(x)$ and $U(x) = B(x)$ absolutely converge on $[0,1)$, and
we have the formal relation
$V(x)=1-\frac{1}{U(x)}$. We bound $u_{2n}$ from below
by Theorem~\ref{thm_effStirl}: for every $n\in\N$,
$$
{\textstyle
u_{2n}=4^{-2n}\frac{(2n)!^2}{n!^4}\ge\frac{\mathrm{e}^{2/(24n+1)}}{\pi n\mathrm{e}^{1/3n}}\ge
\frac{1}{\pi n\mathrm{e}^{1/3}}
}
$$
and $u_{2n}\ge\frac{0.228}{n}$. For $N\in\N$ we set $x=x(N)=1-
N^{-8/9}$ ($\in[0,1)$).
Then for $N\ge 2$,
\begin{eqnarray*}
{\textstyle
\frac{x^{N+1}}{1-x}}&=&{\textstyle\frac{N^{8/9}}{\exp((N+1)\log(1/(1-N^{-8/9})))}}\\
&=&{\textstyle\frac{N^{8/9}}{\exp((N+1)(N^{-8/9}+N^{-2\cdot 8/9}/2+\ds))}\le
\frac{N^{8/9}}{\exp(N^{1/9})}}\,.
\end{eqnarray*}
For $N\ge 3$ we have $x\in 
(\frac{1}{2},1)$ and 
$\log(\frac{1}{1-x})=\frac{8}{9}\log N$.
We use Proposition~\ref{prop_effWeakAbel1} 
with $c=0.228$. Then 
$\frac{8}{9}c\ge0.2$ and $c\log 2\le0.16$. The equality $\frac{a_n}{4^n}
=\sum_{n=0}^N\frac{c_n}{4^n}$ in the proof of part~1 
of Proposition~\ref{prop_firstRel} gives the stated bound that for 
every integer $N\ge 3$,
$$
{\textstyle
1\ge 4^{-N}|W_{\mathrm{2,rec}}(N)|=\sum_{n=0}^N v_n\ge1-\frac{1}{0.2\log N-0.16}-\frac{N^{8/9}}{\exp(N^{1/9})}\,.
}
$$
\eproof

We adapt Proposition~\ref{prop_effWeakAbel1} 
for the bound $u_{2n}\le\frac{c'}{n}$.

\begin{prop}\label{prop_effWeakAbel2}
Suppose that the formal power series
$${\textstyle
U(x)=\sum_{n=0}^{\infty}u_nx^n 
\,\text{ and }\,
V(x)=\sum_{n=0}^{\infty}v_nx^n
}
$$ 
in $\R_{\ge0}[[x]]$  absolutely converge for every
$x\in[0,1)$. Let $u_0=1$, $u_{2n-1}=0$ and $u_{2n}\le\frac{c'}{n}$ for every $n\in\N$ and some constant $c'>0$, and let formally
$V(x)=1-\frac{1}{U(x)}$. Then for every $N\in\N$ and $x\in(0,1)$ with $2N(1-x)\le\frac{1}{3}$ we have
$$
\sum_{n=0}^N v_n\le 1+3N(1-x)-\frac{1}{\big(c'+\frac{1}{\log(6N)}\big)\cdot\log\big(\frac{1}{1-x}\big)}\;.
$$
\end{prop}
\proof
By the assumptions,  $1\le\frac{\log(1/(1-x))}{\log(6N)}$ and $x^{2n}\le x^n$. So for every $x$ as stated we have
$${\textstyle
F_U(x)=\sum_{n=0}^{\infty}u_nx^n\le1+c'\sum_{n=1}^{\infty}\frac{x^{2n}}{n}
\le\big(c'+\frac{1}{\log(6N)}\big)
\log\big(\frac{1}{1-x}\big)\,.
}
$$
For every $N$ in $\N$ and $x=1-y \in(0,1)$ with $y\in(0,1)$ 
so close to $0$ that $2Ny\le\frac{1}{3}$ we therefore have, since $u_n\ge 0$, the stated bound 
\begin{eqnarray*}
{\textstyle\sum_{n=0}^N v_n}&\le&{\textstyle x^{-N}\sum_{n=0}^{\infty}v_n x^n
=\frac{1-1/F_U(x)}{x^N}\le
\frac{1}{x^N}\big(1+
\frac{1}{(c'+
\frac{1}{\log(6N)})\log(1-x)}\big)}\\
&=&\mathrm{e}^{N\log(1/(1-y))}
(\cdots)\le\mathrm{e}^{2Ny}
(\cdots)
\le(1+3Ny)(\cdots)\\
&\le&{\textstyle 1+3N(1-x)-
\frac{1}{(c'+\frac{1}{\log(6N)})\log(\frac{1}{1-x})}\,.}
\end{eqnarray*}
\eproof

From this we deduce an effective lower bound on the difference
$$
1-4^{-N}|W_{\mathrm{2,rec}}(N)|
$$ 
in Theorem~\ref{thm_polya1}
for $d=2$. 

\begin{prop}\label{prop_Po_1_2_lb}
If $N\ge140000$ then
$$
1-4^{-N}|W_{\mathrm{2,rec}}(N)|\ge
(0.84\log N)^{-1}-3N^{-1}\;.
$$
\end{prop}
\proof
Again $v_n=c_n4^{-n}$ and 
$u_n = b_n4^{-n}$, so that $V(x)=C(x)$ and $U(x)=B(x)$. The assumptions of the previous
proposition are satisfied if we get the required upper bound on $u_{2n}$. By Theorem~\ref{thm_effStirl} for every $n\in\N$,
$${\textstyle
u_{2n}=4^{-2n}\frac{(2n)!^2}{n!^4}\le
\frac{\mathrm{e}^{1/12n}}{\pi n\mathrm{e}^{4/(12n+1)}}\le
\frac{\mathrm{e}^{1/12}}{\pi n}
}
$$
and $u_{2n}\le\frac{0.346}{n}$. 
We set
$x=x(N)=1-N^{-2}\in [0,1)$. For $N\ge1.4\cdot 10^5$ one has $2N(1-x)\le
\frac{1}{3}$ and 
$\frac{1}{\log(6N)}\le0.074$. For $c'=0.346$ the previous proposition gives for
$4^{-N}a_N=\sum_{n=0}^N v_n$
the stated bound.
\eproof

We combine both bounds in one estimate.

\begin{cor}\label{cor_poslCOr}
Let $d=2$. There is an effective constant $N_0$ such that for every $N\ge N_0$,
$$
(0.9\log N)^{-1}\le
1-4^{-N}|W_{\mathrm{2,rec}}(N)|
\le(0.1\log N)^{-1}\,.
$$   
\end{cor}

\section{The model 
$\mathrm{W_{\C}}$}\label{sec_weights}

In this section we obtain our main results,  Theorems~\ref{thm_polya1gen1} and 
\ref{thm_polya1gen1apul} extending 
 to $\mathrm{W_{\C}}$
Theorem~\ref{thm_explVal1d3} on $\overline{0}$-recurrence, and
Theorems~\ref{thm_polya1gen2n} 
and \ref{thm_polya1gen2vC}  extending to $\mathrm{W_{\C}}$
Theorem~\ref{thm_explVal2d3} on $\overline{v}$-recurrence. By using 
weights $h\cc\N_2\to\C$ we go 
beyond probability. We employ
absolute convergence of complex series which was recalled in
Section~\ref{sec_intro}. 

We begin with the infinite triangle inequality and then introduce linear
combination, grouping and product
of series. The proof for the
first operation is omitted.

\begin{prop}\label{prop_infTri}
Let $\sum_{x\in X}h(x)$
be an absolutely convergent series. Then the sums satify
$$
{\textstyle
\big|\sum_{x\in X}h(x)
\big|\le
\sum_{x\in X}|h(x)|\,.
}
$$
\end{prop}
\proof
Let $\sum_{x\in X}h(x)$ be as stated. For finite $X$ this is the 
usual triangle inequality in $\C$. Suppose that $X$ is
countable and that the series has the sum $s$ ($\in\C$). Then 
also $\sum_{x\in X}|h(x)|$ absolutely converges and has the sum 
$t$ ($\in\R_{\ge0}$). For any bijection $f\cc\N\to X$, 
\begin{eqnarray*}
|s|&=&\big|\lim_{n\to\infty}
{\textstyle
\sum_{j=1}^n h(f(j))}
\big|=\lim_{n\to\infty}
{\textstyle
\big|\sum_{j=1}^n h(f(j))\big|}\\
&\le&\lim_{n\to\infty}{\textstyle
\sum_{j=1}^n|h(f(j))|=t\,,
}
\end{eqnarray*}
due to the continuity of the function $u\mapsto|u|$, $u\in\C$.
\eproof

If $\al,\be\in\C$ and 
$S_i=\sum_{x\in X}h_i(x)$, $i=1,2$, are two series on  a~set $X$, then the series 
$${\textstyle
\al S_1+\be 
S_2=\sum_{x\in 
X}(\al h_1(x)+\be h_2(x))
}
$$
is the {\em linear combination} of $S_1$ and $S_2$.

\begin{prop}\label{prop_ACser1apul}
If $S_1$ and $S_2$ absolutely converge and have respective sums
$s_1$ and $s_2$, then $\al S_1+\be S_2$ absolutely converges
and has the sum $\al s_1+\be s_2$. 
\end{prop}

For the next two operations we need a~lemma.

\begin{lem}\label{lem_apprSum}
Let $S=\sum_{x\in X}h(x)$ be an absolutely convergent
series with the sum $s$. Then for every $\ep>0$ there is a~finite set 
$Y\sus X$, denoted by $Y(S,\ep)$, such that for every finite set $Z$ with $Y\sus Z\sus X$ we have
$$
{\textstyle
\big|s-\sum_{x\in Z}h(x)\big|\le\ep\,.
}
$$
\end{lem}
\proof
For finite $X$ we take $Y=X$. For countable $X$ we take any bijection
$f\cc\N\to X$. For the given $\ep>0$ we take an $N\in\N$ such that $|\sum_{i=1}^N h(f(i))-
s|\le\frac{\ep}{2}$ and  $\sum_{n>N}|h(f(n))|\le\frac{\ep}{2}$. We
set $Y=f[\,[N]\,]$. Then for every finite set $Z$ with $Y\sus Z\sus X$, 
$$
{\textstyle
|\sum_{x\in Z}h(x)-s|\le
|\sum_{x\in Y}h(x)-s|+
\sum_{x\in Z\setminus Y}|h(x)|\le
\frac{\ep}{2}+\frac{\ep}{2}=\ep
\,.
}
$$
\eproof

A~{\em partition} of a~set $X$ is a~set $P$ of nonempty and
disjoint sets such that $\bigcup P=X$. If $X$ is at most countable, 
then so is $P$ and every set $Z\in P$. If $S=\sum_{x\in X}h(x)$ is 
a~series and $P$ is a~partition of $X$ such that for 
every $Z\in P$ the subseries $\sum_{x\in Z}h(x)$ of $S$ absolutely converges and has the 
sum $s_Z$, then the series
$$
{\textstyle
S_P=\sum_{Z\in P}s_Z
}
$$
is the {\em $P$-grouping} of $S$. 

\begin{prop}\label{prop_grouSer}
Suppose that the series $S=\sum_{x\in X}h(x)$ absolutely
converges, has the sum $s$ and that $P$ is a~partition of $X$. Then
$S_P$ is correctly defined, absolutely
converges, and has the same sum $s$
as $S$. 
\end{prop}
\proof
Let $S$, $s$ and $P$ be as stated. For every $Z\in P$ the series
$S_Z=\sum_{x\in Z}h(x)$ is a~subseries of $S$ and absolutely converges. 
Thus the series $S_P$ is correctly defined.

We show that $S_P$ absolutely converges. Let
$${\textstyle
c=\sup(\{\sum_{x\in Z}|h(x)|\cc\;\text{$Z\sus X$ and is finite}\})\ \ (\in[0,\,+\infty))
}
$$
and let $\{Z_1,Z_2,\ds,Z_n\}\sus P$ be a~finite set. Let $i\in[n]$. We use Lemma~\ref{lem_apprSum} and take the finite set $Z_i'=
Y(S_{Z_i},2^{-i})$ ($\sus Z_i$). We set $Z_0=Z_1'\cup\ds\cup Z_n'$ ($\sus X$); this is a~disjoint union. Then, in the sense of sums,
\begin{eqnarray*}
{\textstyle
\sum_{i=1}^n\big|\sum_{x\in Z_i}h(x)\big|
}&\le&   
{\textstyle
\sum_{i=1}^n\big|\sum_{x\in Z_i}h(x)-\sum_{x\in Z_i'}h(x)\big|\,+
}\\
&+&{\textstyle
\sum_{x\in Z_0}|h(x)|\le\sum_{i=1}^n
2^{-i}+c
}\\
&\le& 1+c\,.
\end{eqnarray*}
Thus $S_P$ absolutely converges. 

Let $t$ be the sum of $S_P$.
We show that $t=s$.  It suffices 
to show that for any given $\ep$ we have $|t-s|\le\ep$. We use 
Lemma~\ref{lem_apprSum} 
and take finite sets $X'\equiv 
Y(S,\frac{\ep}{3})$ ($\sus X$)
and $P'\equiv Y(S_P,\frac{\ep}{3})$ ($\sus P$). Then we take a~finite 
superset of blocks $\{Z_1,\ds,Z_n\}$ of $P'$, that is $P'\sus
\{Z_1,\ds,Z_n\}\sus P$
with $n\in\N$, such that also $X'\sus\bigcup_{i=1}^n Z_i$. We use Lemma~\ref{lem_apprSum} for
every $i\in[n]$ and 
take the finite set $Z_i'= 
Y(S_{Z_i},2^{-i}\frac{\ep}{3})$ ($\sus Z_i$). Finally, for every
$i\in[n]$ we set $Z_i''= Z_i'\cup(X'\cap Z_i)$ ($\sus Z_i$) 
and take the disjoint union $X_0=\bigcup_{i=1}^n Z_i''$ 
($\sus X$). Thus $X_0$ is finite,  $X'\sus X_0$, $Z_i'\sus Z_i''$ and
\begin{eqnarray*}
|s-t|&\le&{\textstyle|s-\sum_{x\in X_0}h(x)|+
\sum_{i=1}^n|\sum_{x\in Z_i''}h(x)-\sum_{x\in Z_i}h(x)|\,+}\\
&+&{\textstyle |\sum_{i=1}^n \sum_{x\in Z_i}h(x)-t|\le\frac{\ep}{3}+\frac{\ep}{3}+\frac{\ep}{3}=\ep\,.
}
\end{eqnarray*}
Hence $s=t$.
\eproof

If $S_i=\sum_{x\in X_i}h_i(x)$, $i=1,2$, are two series on two sets $X_i$, then the 
series 
$${\textstyle
S_1\cdot S_2=\sum_{(x,y)\in 
X_1\times X_2}h_1(x)h_2(y)
}
$$ 
is 
the {\em product} of $S_1$ and $S_2$. 

\begin{prop}\label{prop_proSer}
If $S_1$ and $S_2$ absolutely converge and have respective sums
$s_1$ and $s_2$, then $S_1\cdot S_2$ absolutely converges 
and has the sum $s_1s_2$.
\end{prop}
\proof
We show that $S_1\cdot S_2$ absolutely converges. Suppose that 
$c\ge0$ is such that for $i=1,2$ and any finite sets $Y_i$ with $Y_i\sus 
X_i$ one has
$\sum_{x\in Y_i}|h_i(x)|\le c$.
Let $A\sus X_1\times X_2$ be a~finite set. We 
take any finite sets $Y_i\sus X_i$ such that $A\sus Y_1\times Y_2$ and 
get the desired bound:
$${\textstyle
\sum_{(x,\,y)\in A}|h_1(x)h_2(y)|\le
\sum_{x\in Y_1}|h_1(x)|\cdot
\sum_{y\in Y_2}|h_2(y)|\le c\cdot c=c^2\,.
}
$$

We show that the sum $s$ of $S_1\cdot S_1$ equals $s_1s_2$. 
It suffices to prove that for any given $\ep\in(0,1]$ we have $|s-
s_1s_2|\le\ep$. We use Lemma~\ref{lem_apprSum} and take
finite sets $X'_1=Y(S_1,\frac{\ep}{4(1+|s_2|)})$ ($\sus X_1$), 
$X'_2=Y(S_2,\frac{\ep}{4(1+|s_1|)})$ ($\sus X_2$) and $Z=Y(S_1\cdot 
S_2,\frac{\ep}{4})$ ($\sus X_1\times X_2$). We take finite 
sets $X_1''$ and $X_2''$ such that $X_1'\sus X_1''\sus X_1$, $X_2'\sus 
X_2''\sus X_2$ and $Z\sus X_1''\times X_2''$. Then $|s-
s_1s_2|$ is at most
\begin{eqnarray*}
&&
{\textstyle
\big|s-\sum_{(x,y)\in X_1''\times X_2''}h_1(x)h_2(y)\big|+
\big|\sum_{x\in X_1''}h_1(x)
\sum_{y\in X_2''}h_2(y)-s_1s_2\big|
}\\
&&
{\textstyle
\le\frac{\ep}{4}+|(s_1+t_1)(s_2+t_2)-s_1s_2|,\ |t_1|\le\frac{\ep}{4(1+|s_2|)}
\,\text{ and }\,
|t_2|\le\frac{\ep}{4(1+|s_1|)}\,.
}
\end{eqnarray*}
Hence $|s-s_1s_2|\le
\frac{\ep}{4}+|s_1t_2|+|s_2t_1|+|t_1t_2|
\le\frac{\ep}{4}+\frac{\ep}{4}+
\frac{\ep}{4}+\frac{\ep}{4}=\ep$.
\eproof

As we know, $W(n)$ denotes the set of $1$-walks of length $n$ 
($\in\N_0$) in $K_{\N}$. These are 
the tuples $w=\langle u_0,u_1,\ds,u_n\rangle$, $u_i\in\N$, such that $u_0=1$ and 
$u_{i-1}\ne u_i$ for every $i\in[n]$. Recall that a~weight $h\cc\N_2\to\C$ is light if for every $n\in\N_0$ the
series $\sum_{w\in W(n)}h(w)$ absolutely converges. 

We obtain a~sufficient condition for lightness in terms of values $h(e)$, 
$e\in\N_2$, of weights. We use it later in the follow-up article 
\cite{klaz_hors2}. Let $h\cc\N_2\to\C$ and $n\in\N_0$. By 
$V(h,n)$ ($\sus V(h)$) we denote 
the set of vertices in $G(h)$ lying in distance at most $n$ from the 
vertex $1$.

\begin{defi}[slimness]\label{def_slimWei}
A~weight $h\cc\N_2\to\C$ is slim if for every $n\in\N_0$ there is 
a~constant $c=c(n)\ge0$ such that for every vertex 
$u\in V(h,n)$ and every finite set $X\sus\N\setminus\{u\}$ we have
$${\textstyle
\sum_{v\in X}|h(\{u,\,v\})|\le c\,.
}
$$
\end{defi}
For example, if every vertex in $G(h)$ has finite degree then $h$ is slim.

\begin{prop}\label{prop_slimWei}
Every slim weight is light.
\end{prop}
\proof
Let $h\cc\N_2\to\C$ be a~slim weight. We show by induction on $n\in\N_0$ 
that for every $n$ the series $\sum_{w\in W(n)}h(w)$
absolutely converges. For $n=0$ this is trivial. Let $n>0$ and let
$c\ge0$ be a~constant such that for every finite set $X\sus W(n-1)$ we have 
$\sum_{w\in X}|h(w)|\le c$ and that
for every vertex $u\in V(h,n-1)$ and every 
finite set $X\sus\N\setminus\{u\}$ we have $\sum_{v\in X}|h(\{u,\,v\})|\le c$. 
Let $X\sus W(n)$  
be a~finite set. We may assume that every edge in every walk $w\in X$ 
has nonzero weight. Recall that $\ell(w)$ is the last vertex of a~walk $w$. We 
decompose every walk $w\in X$ as $w=w'\ell(w)$ with $w'\in W(n-1)$. We
denote the set of the walks $w'$ by $Y$, and for any $w'\in Y$ we set $X(w')=\{v\in\N\setminus
\{\ell(w')\}:\;w'v\in X\}$.  
By Propositions~\ref{prop_grouSer} and \ref{prop_ACser1apul}, 
\begin{eqnarray*}
{\textstyle
\sum_{w\in X}|h(w)|}&=&{\textstyle
\sum_{w'\in Y}|h(w')|\cdot\sum_{v\in X(w')}|h(\{\ell(w'),\,v\})|
}\\
&\le&{\textstyle
\sum_{w'\in Y}|h(w')|\cdot c\le c^2\,.
}
\end{eqnarray*}
\eproof

Let $h\cc\N_2\to\C$ be a~weight and $W$ be a~set of walks in 
$K_{\N}$. If the series $\sum_{w\in W}h(w)$ absolutely converges, we denote its 
sum by $h(W)$ ($\in\C$). For $W=\emptyset$ we set $h(W)=0$.
For an extension of  
Theorem~\ref{thm_explVal1d3} on $\overline{0}$-recurrence to $\mathrm{W_{\C}}$ we introduce 
analogs of the counting sequences $(a_n)$, $\ds$, 
$(d_n)$ and generating functions $A(x)$, $\ds$, $D(x)$. Let $h\cc\N_2\to\C$ be 
a~light weight and $n\in\N_0$. We define
three sets of $1$-walks of length $n$ in $K_{\N}$: 
$W_{\mathrm{rec}}(n)$ are the recurrent walks revisiting the 
start $1$, $W_b(n)$ are the walks 
that end at $1$, and $W_{c}(n)$ are the walks that end at $1$ 
and the inner avoids $1$. $W(n)$ is of course the set of all 
$1$-walks of length $n$ in $K_{\N}$.
We define $a_n^h=h(W_{\mathrm{rec}}(n))$, $b_n^h=h(W_b(n))$,
 $c_n^h=h(W_c(n))$ and 
$d_n^h=h(W(n))$. These numbers are correctly defined because $h$ is 
light and all four series are subseries of $\sum_{w\in 
W(n)}h(w)$. We have $a_0^h=c_0^h=0$ and 
$b_0^h=d_0^h=1$. We define the generating functions 
\begin{eqnarray*}
&&{\textstyle
A_h(x)=\sum_{n=0}^{\infty}a_n^hx^n,\ 
B_h(x)=\sum_{n=0}^{\infty}
b_n^hx^n,\ C_h(x)=\sum_{n=0}^{\infty}
c_n^hx^n\,\text{ and }}\\ 
&&{\textstyle
D_h(x)=\sum_{n=0}^{\infty}
d_n^hx^n\,.}
\end{eqnarray*}
All four are in $\C[[x]]$. 
The quantities of primary interest are $a_n^h$ and $A_h(x)$. We 
generalize Proposition~\ref{prop_firstRel}. 

\begin{prop}\label{prop_12rel}
Let $h\cc\N_2\to\C$ be a~light weight and the generating functions $A_h(x)$, $B_h(x)$, 
$C_h(x)$ and $D_h(x)$ be as above. The following formal relations hold in $\C[[x]]$.
\begin{enumerate}
\item $A_h(x)=C_h(x)D_h(x)$. 
\item $B_h(x)=\frac{1}{1-C_h(x)}=
\sum_{j=0}^{\infty}C_h^j(x)$, equivalently $C_h(x)=1-
\frac{1}{B_h(x)}$.
\end{enumerate}
\end{prop}
\proof
1. Let $n\in\N_0$. We show that
$a_n^h=\sum_{j=0}^n
c_j^h\cdot d_{n-j}^h$. For $n\le1$
it is trivial, then $a_n^h=c_n^h=0$. For $n\ge2$ we have
$W_{\mathrm{rec}}(n)\ne\emptyset$ and
$W_{\mathrm{rec}}(n)$ is in fact countable. We use the bijection
$${\textstyle
F\cc W_{\mathrm{rec}}(n)\to
\bigcup_{j=0}^n W_c(j)\times W(n-j)
}
$$
that is defined as follows. For 
$w=\langle u_0,\ds,u_j,\ds,u_n\rangle$ in $W_{\mathrm{rec}}(n)$, where 
$u_0=u_j=1$, $j>0$ and $u_i\ne1$ for $i\in[j-1]$, we set 
$F(w)=(w_1,w_2)$, where
$w_1=\langle u_0,u_1,\ds,u_j\rangle$ and $w_2=\langle u_j,u_{j+1},\ds,u_n\rangle$ are uniquely determined. 
The map $F$ is weight-preserving, or  WP, in the sense that for 
every weight $h\cc\N_2\to\C$ and 
value $F(w)=(w_1,w_2)$ we have $h(w)=h(w_1) 
h(w_2)$. We extend $h$ to pairs of walks by setting 
$h((w,w'))=h(w)h(w')$. Let $h\cc\N_2\to\C$ be a~light weight. Using the bijection $F$ and 
the above propositions we indeed get 
\begin{eqnarray*}
a_n^h&=&h(W_{\mathrm{rec}}(n))
\stackrel{\mathrm{WP}}{=}h(F[W_{\mathrm{rec}}(n)])\\
&\stackrel{\text{Prop.~\ref{prop_grouSer}}}{=}&{\textstyle\sum_{j=0}^n h\big(W_c(j)\times W(n-j)\big)}\\
&\stackrel{\text{Prop.~\ref{prop_proSer}}}{=}&{\textstyle\sum_{j=0}^nh(W_c(j))\cdot h(W(n-j))
=\sum_{j=0}^nc_j^h\cdot d_{n-j}^h}\,.
\end{eqnarray*} 

2. Let $n\in\N_0$. We show that 
$$
{\textstyle
b_n^h=\sum_{j=0}^{\infty}\sum_{\substack{n_1,\,\ds,\,n_j\in\N_0\\n_1+
\ds+n_j=n}}c_{n_1}^h c_{n_2}^h\ds c_{n_j}^h\,,
}
$$
where the empty sum is set to $0$,  
and the empty product to $1$. For 
every $n$ the double sum has only 
finitely many summands. This equality again holds trivially for 
$n\le 1$. If $n=0$ then 
$b_0^h=1$ and this $1$ is 
provided by the summand on the right-hand side with $j=0$ that is defined
for $n=0$ as $1$, and for $n>0$ as $0$. If $n=1$
then $b_1^h=c_1^h=0$. For $n\ge2$ we have $W_b(n)\ne\emptyset$ and 
$W_b(n)$ is countable. We use the bijection
$${\textstyle
F\cc W_b(n)\to
\bigcup_{j=1}^{\infty}\bigcup_{\substack{n_1,\,\ds,\,n_j\in\N\\n_1+\ds+n_j=n}}W_c(n_1)\times\cdots\times W_c(n_j)
}
$$
that is defined as follows. For 
$$
w=\langle u_0,\,\ds,\,u_{m_1},\,\ds,\,\ds,\,u_{m_{j-1}},\,\ds,\,u_n\rangle\in W_b(n)\,,
$$ 
where $0=m_0<m_1<\ds<m_j=n$ and $u_k=1$ iff $k=0$ or $k=m_i$ for some $i\in[j]$, we set
$F(w)=\langle w_1,w_2,\ds,w_j\rangle$ where  
$$
w_i=\langle u_{m_{i-1}},\,u_{m_{i-1}+1},\,\ds,\,u_{m_i}\rangle,\ i\in[j] 
$$
(thus $n_i=m_i-m_{i-1}$). 
Again $F$ is weight-preserving, or WP. We formally 
extend any weight $h\cc\N_2\to\C$ to tuples of walks by 
$${\textstyle
h(\langle w_1,\,w_2,\,\ds,\,w_j\rangle)=
\prod_{i=1}^j h(w_i)\,. 
}
$$
Let $h$ be a~light weight. Using the bijection $F$ and the above 
propositions we indeed have 
\begin{eqnarray*}
b_n^h&=&h(W_b(n))\stackrel{\mathrm{WP}}{=}h(F[W_b(n)])\\
&\stackrel{\text{Prop.~\ref{prop_grouSer}}}{=}&{\textstyle\sum_{j=1}^{\infty}\sum_{\substack{n_1,\,\ds,\,n_j\in\N_0\\n_1+
\ds+n_j=n}}h\big(W_c(n_1)\times\cdots\times W_c(n_j)\big)}\\
&\stackrel{\text{Prop.~\ref{prop_proSer}}}{=}&{\textstyle\sum_{j=1}^{\infty}\sum_{\cdots}h(W_c(n_1))\ds h(W_c(n_j))
=\sum_{j=1}^{\infty}\sum_{\cdots}c_{n_1}^h\ds c_{n_j}^h}\,.
\end{eqnarray*}
\eproof

\noindent
This proof uses 
our {\em semi-formal approach}, or SFA, to generating functions. SFA extends the symbolic method in 
(enumerative) combinatorics, or SMC, 
\cite{albe_al,comt,flaj_sedg,goul_jack,symb} from finite to countable 
sets. For example, the product formula for generating functions in SMC is based on the identities 
$$
|A\times B|=|A|\cdot|B|\,\text{ and }\,|C\cup D|=|C|+|D|
$$
for cardinalities of finite sets $A$, $B$, $C$ and $D$, where $C\cap 
D=\emptyset$. For a~formal statement
of the product formula
see, for example, 
\cite[2.2.14]{goul_jack}. SFA extends it to identities like the 
above $A_h(x)=C_h(x)D_h(x)$. The product formula for generating functions in SFA is 
based on Propositions~\ref{prop_proSer} and
\ref{prop_grouSer} which deal with at most countable sets. We plan to 
discuss SFA in more detail in \cite{klaz1}.

We call a~light weight $h\cc\N_2\to\C$  
{\em convex} if for every vertex $u\in V(h)$ we have the sum
$${\textstyle
\sum_{v\in V(h)\setminus\{u\}}h(\{u,\,v\})=\sum_{v\in\N\setminus\{u\}}h(\{u,\,v\})=1\,.
}
$$
We show that these sums exist. If they are all $\ge 1$, we call $h$ {\em superconvex}.

\begin{prop}\label{prop_convFromLigh}
Let $h\cc\N_2\to\C$ be a~light weight and $u\in V(h)$. 
Then the series $\sum_{v\in\N\setminus\{u\}}
h(\{u,v\})$ absolutely converges.   
\end{prop}
\proof
Let $h$ and $u$ be as stated. There is a~$1$-walk $w_0=\langle 
u_0,u_1,\ds,u_n\rangle$ in $G(h)$ such that $u_n=u$. We set 
$${\textstyle
c=h(w_0)=\prod_{i=1}^n h(\{u_{i-1},u_i\})\ \  (\ne0) 
}
$$
and denote by $W$ the set of $1$-walks in $K_{\N}$ of the form 
$$
\langle u_0,\,u_1,\,\ds,\,u_n,\,u_{n+1}\rangle\,. 
$$
Thus $u_0$, $u_1$, $\ds$, $u_n$ are fixed and $u_{n+1}$ runs in 
$\N\setminus\{u_n\}=\N\setminus\{u\}$. Then $\sum_{w\in W}h(w)$ absolutely 
converges because it is 
a~subseries of the absolutely convergent series $\sum_{w\in W(n+1)}h(w)$. By 
Proposition~\ref{prop_ACser1apul} the series 
$$
{\textstyle
\sum_{v\in\N\setminus\{u\}}
h(\{u,v\})=\frac{1}{c}\sum_{w\in W}h(w)
}
$$
absolutely converges. 
\eproof

\noindent
Convexity is a~remnant of the probabilistic approach. Using it we can generalize the identity $D(x)=\frac{1}{1-x}$.

\begin{prop}\label{prop_onD}
For every convex light weight $h\cc\N_2\to\C$ we have
$${\textstyle
D_h(x)=\sum_{n=0}^{\infty}h(W(n))x^n=\sum_{n=0}^{\infty}x^n
=\frac{1}{1-x}
}\,.
$$ 
\end{prop}
\proof
Let $h$ be a~convex light weight. We prove by induction on $n\in\N_0$ that 
$h(W(n))=d_n^h=1$. 
For $n=0$ this holds as $W(0)=\{(1)\}$. Let $n>0$. Recall that $\ell(w)$ denotes 
the last vertex of a~walk. By
Propositions~\ref{prop_grouSer} and \ref{prop_ACser1apul}, the convexity and induction we get that indeed
\begin{eqnarray*}
h(W(n))&=&
{\textstyle \sum_{w\in W(n)}h(w)
}\\
&=&{\textstyle\sum_{\substack{w'\in W(n-1)\\\ell(w')\in V(h)}}h(w')\sum_{v\in\N\setminus\{\ell(w')\}}h(\{\ell(w'),\,v\})}\\
&=&{\textstyle\sum_{w'\in W(n-1)}h(w')\cdot1=1\,.
}
\end{eqnarray*}
We can add and remove the condition $\ell(w')\in V(h)$ without affecting
the sum because the walks $w'$ not satisfying it have zero weights. 
\eproof

We generalize Definition~\ref{def_UzR}.

\begin{defi}[$U(R)$ 2]\label{def_UzRvC}
Let $R>0$ be a~real number and  $U(x)=\sum_{n=0}^{\infty}u_nx^n$ be 
in $\C[[x]]$. We define the complex number
$${\textstyle
U(R):=\lim_{n\to\infty}\sum_{j=0}^nu_jR^j\ \ 
(=:\sum_{n=0}^{\infty}u_jR^j
})\,,
$$
if this limit exists.
\end{defi}
This is slightly more 
restrictive than 
Definition~\ref{def_UzR}, we do not allow the value $U(R)=+\infty$.

We obtain a~complex analog of Proposition~\ref{prop_tauber1}.

\begin{prop}\label{prop_tauber1vC}
Suppose that $U(x)=\sum_{n=0}^{\infty}u_n x^n$ and $V(x)=\sum_{n=0}^{\infty}v_n x^n$ 
are in $\C[[x]]$, that the sum $V(1)$ exists and that formally $U(x)=V(x)\frac{1}{1-x}$.  
Then
$$
\lim_{n\to\infty}u_n=
\lim_{n\to\infty}(v_0+v_1+\ds+v_n)=V(1)\ \ (\in\C)\,.
$$
\end{prop}
\proof
Again, $u_n=v_0+v_1+\ds+v_n$ for every $n\in\N_0$.
\eproof

We make use of the 
next classical theorem, see \cite[Chapter II.7]{tene}. For 
a~real number $R>0$ we denote by 
$$
K_R=\{x\in\C:\;|x|<R\}
$$ 
the open disc in $\C$ with 
radius $R$ and center at $0$. It is well known that if $U(x)$ is in 
$\C[[x]]$ and $U(R)$ exists, then 
the series $U(x)$ absolutely converges for every $x\in K_R$.

\begin{thm}[Abel's]\label{thm_abel}
Let $R>0$ and $U(x)=\sum_{n=0}^{\infty}u_nx^n$ 
be in $\C[[x]]$ such that the sum $U(R)$ exists. Then the sum 
$F_U\cc K_R\to\C$ of $U(x)$ has the limit 
$$
\lim_{\substack{x\to R\\x\in[0,\,R)}}F_U(x)=U(R)\ \ (\in\C)\,.
$$
\end{thm}

\begin{cor}\label{cor_tauber1apulvC}
Suppose that $U(x)=\sum_{n=0}^{\infty}u_n x^n$, $V(x)=\sum_{n=0}^{\infty}v_n x^n$ and $W(x)=\sum_{n=0}^{\infty}w_n x^n$ are in $\C[[x]]$  such that formally $U(x)=V(x)W(x)$. If the sums $U(1)$,
$V(1)$ and $W(1)$ exist, then 
$$
U(1)=V(1)W(1)\ \ (\in\C)\,.
$$
\end{cor}
\proof
Propositions~\ref{prop_grouSer} and \ref{prop_proSer} 
imply that $F_U(x)=F_V(x)F_W(x)$ for every $x\in K_1$. Thus by 
Theorem~\ref{thm_abel} and by properties of limits of functions,
\begin{eqnarray*}
U(1)&=&\lim_{\substack{x\to1
\\x\in[0,\,1)}}F_U(x)=
\lim_{\substack{x\to1
\\x\in[0,\,1)}}(F_V(x)F_W(x))\\
&=&\lim_{\substack{x\to1
\\x\in[0,\,1)}}F_V(x)
\lim_{\substack{x\to1
\\x\in[0,\,1)}}F_W(x)=V(1)W(1)\,.
\end{eqnarray*}
\eproof

\noindent
In contrast to the nonnegative real case mentioned in the
remark after Corollary~\ref{cor_UVW1}, 
in the complex case one cannot avoid
Abel's theorem and get the equality 
$U(1)=V(1)W(1)$ only by manipulating 
series. The reason is that
the involved series need not be absolutely convergent. 

\begin{cor}\label{cor_UV1vC}
Suppose that $U(x)$ and $V(x)$ are in $\C[[x]]$ such that formally
$V(x)=1-\frac{1}{U(x)}$. If the sums $U(1)$ and $V(1)$ exist, then $U(1)\ne0$ and 
$${\textstyle
V(1)=1-\frac{1}{U(1)}\ \ (\in\C)\,. 
}
$$
\end{cor}
\proof
We have formally $V(x)U(x)=U(x)-1$. We argue, using 
Theorem~\ref{thm_abel}, as in the previous proof.
\eproof

\begin{cor}\label{cor_UVW1vC}
Suppose that $U(x)$, $V(x)$ and $W(x)$ are formal power series in $\C[[x]]$  
such that formally $W(x)=V(x)+W(x)\cdot U(x)^2$. If
the sums $U(1)$, $V(1)$ and $W(1)$ exist, then
$$
W(1)=V(1)+W(1)\cdot U(1)^2\ \ (\in\C)\,.
$$
\end{cor}
\proof
We argue as in the two previous proofs.
\eproof

We arrive at the main result in this section which is the first extension 
of Theorem~\ref{thm_explVal1d3} on $\overline{0}$-recurrence to 
$\mathrm{W_{\C}}$. It concerns general weights $h\cc\N_2\to\C$. We denote by 
$$
|h|\cc\N_2\to\R_{\ge0}
$$ 
the weight 
$|h|(\{m,n\})=|h(\{m,n\})|$.
We divide the theorem in two cases according to the existence of
the sum $D_h(1)$.

\begin{thm}[general $1$-recurrence]\label{thm_polya1gen1}
Let $h\cc\N_2\to\C$ be a~light weight. In the above
notation the following holds.
\begin{enumerate}
\item (In this case $D_h(1)$ exists.) Suppose that the sums 
$D_h(1)$, $C_h(1)$, $B_h(1)$ and $A_h(1)$ exist\,---\,a~sufficient 
condition for it is $D_{|h|}(1)<+\infty$. Then $B_h(1)\ne0$ and
$$
A_h(1)=
{\textstyle
\sum_{n=0}^{\infty}\sum_{w\in W_{\mathrm{rec}}(n)}h(w)=
\big(1-\frac{1}{B_h(1)}\big)D_h(1)\ \ (\in\C)\,.
}
$$
\item 
(In this case $D_h(1)$ does not exist.) Suppose that $h\cc\N_2\to\R_{\ge0}$ and $D_h(1)=+\infty$. Then $c_2^h>0$
and $a_{n+2}^h\ge c_2^hd_n^h$ for every $n\in\N_0$. Hence $A_h(1)=+\infty$.
\end{enumerate} 
\end{thm}
\proof
1. If $D_{|h|}(1)<+\infty$ then the four series with the sums $D_h(1)$, 
$C_h(1)$, $B_h(1)$ and $A_h(1)$ absolutely converge. This follows from 
the fact that for every $n\in\N_0$,
$$
|d_n^h|,\,|c_n^h|,\,|b_n^h|,\,|a_n^h|
\le d_n^{|h|}\,. 
$$
These inequalities follow from Proposition~\ref{prop_infTri} and the 
fact that for every $n$ in $\N_0$ the complex series $\sum_{w\in W(n)}h(w)$, 
$\sum_{w\in W_c(n)}h(w)$, $\sum_{w\in W_b(n)}h(w)$ and $\sum_{w\in 
W_a(n)}h(w)$ are subseries of 
$\sum_{w\in W(n)}h(w)$. 

Suppose that $D_h(1)$, 
$C_h(1)$, $B_h(1)$ and $A_h(1)$ exist. 
By part~1 of 
Proposition~\ref{prop_12rel} we have formally $A_h(x)=C_h(x)D_h(x)$.
By Corollary~\ref{cor_tauber1apulvC}, 
$A_h(1)=C_h(1)D_h(1)$ ($\in\C$).
By part~2 of Proposition~\ref{prop_12rel} and by
Corollary~\ref{cor_UV1vC}, $B_{h}(1)\ne0$ and $C_h(1)=1-\frac{1}{B_h(1)}$. The stated formula follows. 

2. Suppose that $h$ is nonnegative and $D_h(1)=+\infty$. Thus $d_n^h>0$
for some $n>0$ and there must be a~vertex $u\in\N\setminus\{1\}$ 
with $h(\{1,u\})>0$. The walk $\langle 1,u,1\rangle$ then shows that $c_2^h>0$.
Since by part~1 of 
Proposition~\ref{prop_12rel} we have $a_n^h=\sum_{j=0}^n c_j^hd_{n-
j}^h$ for every $n\in\N_0$, the stated inequality follows. Hence
$A_h(1)=+\infty$.
\eproof

\noindent
We remark that we cannot completely characterize either case of the
theorem. The same holds for the remaining three main theorems in
this section.

The first of them is another extension of Theorem~\ref{thm_explVal1d3} on $\overline{0}$-recurrence to
$\mathrm{W_{\C}}$, now for convex weights. We divide the theorem in
two cases according to the existence of 
$B_h(1)$.

\begin{thm}[convex $1$-recurrence]\label{thm_polya1gen1apul}
Let $h\cc\N_2\to\C$ be any convex light weight. In the above notation the following holds.
\begin{enumerate}
\item (In this case $B_h(1)$ exists.) Suppose that the sums $B_h(1)$ and $C_h(1)$ exist\,---\,a~sufficient condition for it is $B_{|h|}(1)<+\infty$. Then $B_{h}(1)\ne0$ and
$$
\lim_{n\to\infty}a_n^h=
\lim_{n\to\infty}h(W_{\mathrm{rec}}(n))=
{\textstyle
1-\frac{1}{B_h(1)}\ \ (\in\C)\,.
}
$$
\item 
(In this case $B_h(1)$ does not exist.) Suppose that
$h\cc\N_2\to\R_{\ge0}$
and $B_h(1)=+\infty$. Then 
$$
\lim_{n\to\infty}a_n^h=
\lim_{n\to\infty}
{\textstyle
\sum_{w\in W_{\mathrm{rec}}(n)}h(w)
}
=1\,.
$$
\end{enumerate} 
\end{thm}
\proof
1. If $B_{|h|}(1)<+\infty$ then we
use the bounds $|b_n^h|,|c_n^h|\le b_n^{|h|}$ and proceed as in the 
proof of part~1 of the previous theorem.

Suppose that $B_h(1)$ and $C_h(1)$ exist.
By part~1 of Proposition~\ref{prop_12rel}, 
Proposition~\ref{prop_onD} and by Proposition~\ref{prop_tauber1vC}, 
$\lim_{n\to\infty}a_n^h=C_h(1)$. By part~2 of Proposition~\ref{prop_12rel} 
and Corollary~\ref{cor_UV1vC}, $B_h(1)\ne0$ and $C_h(1)=1-
\frac{1}{B_h(1)}$. The stated formula follows.

2. Suppose that $h$ is nonnegative and $B_h(1)=+\infty$. By 
Proposition~\ref{prop_onD}, $D_h(x)=\frac{1}{1-x}$. Thus all coefficients $a_n^h$, $b_n^h$ and $c_n^h$ lie in $[0,1]$ 
and formal power series $A_h(x)$, $B_h(x)$ and $C_h(x)$ converge on 
$[0,1)$. By part~1 of Proposition~\ref{prop_12rel}, 
formally $A_h(x)=
C_h(x)\frac{1}{1-x}$. By Proposition~\ref{prop_tauber1}, 
$\lim_{n\to\infty}a_n^h=C_h(1)$.
By part~2 of Proposition~\ref{prop_12rel} and Corollary~\ref{cor_UV1},
$C_h(1)=1$. The value $1$ of the limit follows.
\eproof

\noindent
Unlike in the forthcoming 
Theorems~\ref{thm_polya1gen2n} and \ref{thm_polya1gen2vC} on 
$v$-recurrence, we do 
not restrict weights by
vertex-transitivity. 
This is interesting since we generalized the case of
vertex-transitive graphs $G_d$. 
Theorems~\ref{thm_polya1gen1} and \ref{thm_polya1gen1apul} generalize 
Theorem~\ref{thm_explVal1d3}.

We generalize P\'olya's claim that ``The 
probability obviously grows with $n$.'' Let $n\in\N_0$, $v\in\N$ and 
$w=\langle u_0,u_1\ds,u_n\rangle$. We define
$$
W(v,\,n)=\{w\in W(n):\;\exists\,i\in[n]\,(u_i=v)\}\;.
$$
It is the set of $v$-recurrent $1$-walks of length $n$ in $K_{\N}$. 

\begin{prop}\label{prop_mono}
If $h\cc\N_2\to\R_{\ge0}$ is a~superconvex light weight, then
$$
h(W(v,0))\le h(W(v,1))\le h(W(v,2))\le\cdots
$$
for every $v\in\N$.
\end{prop}
\proof
Let $v\in\N$ and $n\in\N_0$. Recall that $\ell(w)$ is the last vertex of $w$. By Propositions~\ref{prop_ACser1apul} and \ref{prop_grouSer}, 
\begin{eqnarray*}
h(W(v,\,n))&=&{\textstyle\sum_{w\in W(v,\,n)}h(w)=
\sum_{\substack{w\in W(v,\,n)\\
\ell(w)\in V(h)}}h(w)\cdot1}\\
&\le&{\textstyle\sum_{\substack{w\in W(v,\,n)\\\ell(w)\in 
V(h)}}h(w)\cdot\sum_{u\in\N\setminus\{\ell(w)\}}h(\{\ell(w),\,u\})}\\
&=&{\textstyle\sum_{w'\in W(v,\,n+1)}h(w')=h(W(v,\,n+1))\,.}
\end{eqnarray*}
In the second line we used the superconvexity of $h$. We can add and remove the condition $\ell(w)\in V(h)$ because walks not satisfying it have zero weight.
\eproof

We introduce quantities needed to generalize $\overline{v}$-recurrence. Let 
$h\cc\N_2\to\C$ be 
a~light weight, $v\in\N\setminus\{1\}$ and $n\in\N_0$. We define 
sequences of weights of sets of
$1$-walks of length $n$ in $K_{\N}$. Namely,  $(a_n^h)'=h(W(v,n))$, 
$b_n^h=h(W_b(n))$ is as before, $(b_n^h)'=h(W_{b'}(n))$ is the 
weight of walks ending at $1$ and avoiding $v$, $(c_n^h)'=
h(W_{c'}(n))$ is 
the weight of walks ending at $v$ and with the inner avoiding $v$, $(c_n^h)''=
h(W_{c''}(n))$ is the weight of walks ending at $v$ and with the inner avoiding both $1$ and $v$, and 
$d_n^h=h(W(n))$ is as before. 
We have $(a_0^h)'=(c_0^h)'=(c_0^h)''=0$ and 
$b_0^h=(b_0^h)'=d_0^h=1$. Besides the previous generating functions $B_h(x)$ and 
$D_h(x)$ we consider 
\begin{eqnarray*}
&&{\textstyle
A_{0,h}(x)=
\sum_{n=0}^{\infty}(a_n^h)'x^n,\ 
B_{0,h}(x)=\sum_{n=0}^{\infty}
(b_n^h)'x^n,}\\
&&{\textstyle C_{0,h}(x)=
\sum_{n=0}^{\infty}(c_n^h)'x^n \,\text{ and }\,C_{1,h}(x)=
\sum_{n=0}^{\infty}(c_n^h)''x^n\,.
}
\end{eqnarray*} 
All are in $\C[[x]]$. The dependence on $v$ is left 
implicit in the notation.

Let $v\in\N$. We say that a~weight $h\cc\N_2\to\C$ is {\em $v$-transitive} if there is a~bijection 
$f\cc\N\to\N$ such that $f(1)=v$ and for every edge 
$e\in\N_2$ we have 
$h(e)=h(f[e])$. 
We generalize 
Proposition~\ref{prop_firstRel2}.

\begin{prop}\label{prop_firstRel2h}
Let $v\in\N\setminus\{1\}$,  $h\cc\N_2\to\C$ be 
a~$v$-transitive
light weight and let the generating functions
$A_{0,h}(x)$, $B_h(x)$, $B_{0,h}(x)$, $C_{0,h}(x)$, $C_{1,h}(x)$ 
and $D_h(x)$ be as above. The following formal relations hold. 
\begin{enumerate}
    \item $A_{0,h}(x)=
    C_{0,h}(x)D_h(x)$.
    \item $B_h(x)=
    B_{0,h}(x)+
    C_{0,h}(x)^2B_h(x)$,
    equivalently $B_h(x)=
    \frac{B_{0,h}(x)}{1-C_{0,h}(x)^2}$.
    \item $C_{0,h}(x)=
    B_{0,h}(x)C_{1,h}(x)$. 
\end{enumerate}
\end{prop}
\proof
We give less details than in the proof of 
Proposition~\ref{prop_12rel}. 
We argue on the high level of SFA, the semi-formal approach to generating functions. 

1. Any walk $w$ in $\bigcup_{n=0}^{\infty}W(v,n)$ splits at the 
first visit of $v$ in two walks as $w=w_1w_2$, where $w_1$ is 
weight-counted in $C_{0,h}(x)$ and $w_2$ is an arbitrary 
$v$-walk. By the 
$v$-transitivity of $h$, 
the walks $w_2$ are weight-counted by $D_h(x)$. The first relation follows.

2. It suffices to prove the first equality. $B_h(x)$ weight-counts 
$1$-walks $w$ in $K_{\N}$ ending at 
$1$. Those avoiding $v$ are 
weight-counted 
by $B_{0,h}(x)$. Those visiting $v$ split uniquely as $w=w_1w_2w_3$, where the walk $w_1$ starts at $1$, ends at
$v$ and the inner avoids $v$, the walk $w_2$ starts and ends at $v$, and the walk $w_3$ 
starts at $v$, ends at $1$ and the inner avoids $v$. By reversing the
walk $w_3$ we see that both $w_1$ and $w_3$ are
weight-counted by $C_{0,h}(x)$. By the $v$-transitivity of $h$, the 
middle walks $w_2$ are 
weight-counted by $B_h(x)$. The second relation follows.

3. We consider walks $w$ weight-counted by $C_{0,h}(x)$. They 
start at $1$, end at $v$, and the inner avoids $v$. They
uniquely split at the last visit of $1$ in two walks as 
$w=w_1w_2$, 
where $w_1$ and $w_2$ is
weight-counted by $B_{0,h}(x)$
and $C_{1,h}(x)$, respectively. The third relation follows. 
\eproof

The third main result in this section is the first extension of Theorem~\ref{thm_explVal2d3} on
$\overline{v}$-recurrence to $\mathrm{W_{\C}}$. It concerns
general weights $h:\N_2\to\C$. Recall that for any weight 
$h\cc\N_2\to\C$ we denote by $V(h)$ 
($\sus\N$) the vertex set of the component of $1$ in the 
subgraph of $K_{\N}$ formed by edges with nonzero weight. For 
$x\in\C$ we define 
$$
\mathrm{sqrt}(x)=
\{y\in\C\cc\;y^2=x\}\,. 
$$
Always 
$|\mathrm{sqrt}(x)|=2$, except for $\mathrm{sqrt}(0)=\{0\}$. 
For two numbers
$x,z\in\C$ we define 
$$
z\cdot\mathrm{sqrt}(x)=
\{zy\cc\;y\in\mathrm{sqrt}(x)\}\,.
$$
As before we divide the theorem in two cases according to the 
existence of the sum $D_h(1)$, but now we also have the singular case 
when $v\not\in V(h)$.

\begin{thm}[general $v$-recurrence]\label{thm_polya1gen2n}
Suppose that $v\in\N\setminus\{1\}$ and that $h\cc\N_2\to\C$ is a~$v$-transitive 
light weight. In the above notation the following holds.   \begin{enumerate}
\item If $v\not\in V(h)$ then $(a_n^h)'=h(W(v,n))=0$ for every $n\in\N_0$.
\item (In this case $D_h(1)$ exists.) Suppose that the sums $D_h(1)$,
$C_{0,h}(1)$, $B_{0,h}(1)$, $B_h(1)$ and $A_{0,h}(1)$ exist\,---\,a~sufficient condition for it is 
$D_{|h|}(1)<+\infty$. Then $B_h(1)\ne0$ and
$$
A_{0,h}(1)={\textstyle
\sum_{n=0}^{\infty}\sum_{w\in W(v,\,n)}h(w)}
\in
{\textstyle D_h(1)\cdot\mathrm{sqrt}\big(1-\frac{B_{0,\,h}(1)}{B_h(1)}}\big)\ \ 
(\sus\C)\,.
$$
\item (In this case $D_h(1)$ does not exist and $v\in V(h)$.) Let $h\cc\N_2\to\R_{\ge0}$, $D_h(1)=+\infty$ and $v\in V(h)$. 
Then for an $m\in\N$ we have $(c_m^h)'>0$ and $(a_{n+m}^h)'\ge
(c_m^h)'d_n^h$ for every $n\in\N_0$. Hence $A_{0,h}(1)=+\infty$.
\end{enumerate} 
\end{thm}
\proof
Let $v$ and $h$ be as stated.

1. Every walk in $K_{\N}$ joining $1$ and $v$ has weight $0$,
which implies that $(a_n^h)'=0$ for every $n\in\N_0$.

2. For $D_{|h|}(1)<+\infty$ we proceed as in the proof of part~1 of
Theorem~\ref{thm_polya1gen1}. We have bounds 
$$
|d_n^h|,\,|(c_n^h)'|,\,|(b_n^h)'|,\,
|b_n^h|,\,|(a_n^h)'|\le d_n^{|h|}\,.
$$

Suppose that $D_h(1)$, $C_{0,h}(1)$, $B_{0,h}(1)$, $B_h(1)$ and $A_{0,h}(1)$ exist. By part~2 of 
Proposition~\ref{prop_firstRel2h}
and by Corollary~\ref{cor_UVW1vC}, 
$B_h(1)=B_{0,h}(1)
+C_{0,h}(1)^2B_h(1)$. 
From the proof of part~1 of Theorem~\ref{thm_polya1gen1} we 
know that $B_h(1)\ne0$. Hence $C_{0,h}(1)\in\mathrm{sqrt}\big(1-
\frac{B_{0,h}(1)}{B_h(1)}\big)$. By part~1 of Proposition~\ref{prop_firstRel2h},
formally $A_{0,h}(x)=
C_{0,h}(x)D_h(x)$.
By Corollary~\ref{cor_tauber1apulvC}, 
$A_{0,h}(1)=C_{0,h}(1)D_h(1)$ and the stated formula follows.

3. Suppose that $h$ is nonnegative, $D_h(1)=+\infty$ and $v\in V(h)$. 
The shortest walk of length $m$ joining $1$ and $v$ in $G(h)$ shows that $(c_m^h)'>0$. Since by part~1 of 
Proposition~\ref{prop_firstRel2h} we have $(a_n^h)'=\sum_{j=0}^n 
(c_j^h)'d_{n-j}^h$ for every $n\in\N_0$, the stated inequality 
follows. Hence $A_{0,h}(1)=+\infty$.
\eproof

The fourth and last main result in this section is
the second extension of Theorem~\ref{thm_explVal2d3} 
on $\overline{v}$-recurrence to
$\mathrm{W_{\C}}$, now for convex weights. We again divide the theorem 
in two parts according to the existence of $B_h(1)$, and  have the 
singular case when $v\not\in V(h)$.

\begin{thm}[convex $v$-recurrence]\label{thm_polya1gen2vC}
Suppose that $v\in\N\setminus\{1\}$ and that $h\cc\N_2\to\C$ is a~convex 
$v$-transitive light weight. In the above notation the 
following holds.
\begin{enumerate}
\item If $v\not\in V(h)$ then $(a_n^h)'=h(W(v,n))=0$ for every $n\in\N_0$.
\item (In this case $B_h(1)$ exists.)
Suppose that the sums $B_h(1)$, $B_{0,h}(1)$ and $C_{0,h}(1)$ exist\,---\,a~sufficient condition for it is $D_{|h|}(1)<+\infty$. Then $B_h(1)\ne0$ and
$$
\lim_{n\to\infty}(a_n^h)'=
\lim_{n\to\infty}h(W(v,\,n))
\in
{\textstyle\mathrm{sqrt}\big(1-\frac{B_{0,\,h}(1)}{B_h(1)}}\big)\ \ (\sus\C)
\,.
$$
\item (In this case $B_h(1)$ does not exist and $v\in V(h)$.) Let $h\cc\N_2\to\R_{\ge0}$, $B_h(1)=+\infty$
and $v\in V(h)$. Then
$$
\lim_{n\to\infty}(a_n^h)'=
\lim_{n\to\infty}
{\textstyle
\sum_{w\in W(v,\,n)}h(w)
}
=1\,.
$$
\end{enumerate} 
\end{thm}
\proof 
Let $v$ and $h$ be as stated.

1. Again $(a_n^h)'=0$ for every $n\in\N_0$ because $v$ is unreachable
from $1$ in $G(h)$.

2. If $D_{|h|}(1)<+\infty$ then we proceed as in the proof of part~1 of
Theorem~\ref{thm_polya1gen1}. Now we have the bounds 
$$
|b_n^h|,\,|(b_n^h)'|,\,|(c_n^h)'|\le d_n^{|h|}\,.
$$ 

Suppose that $B_h(1)$, $B_{0,h}(1)$ and $C_{0,h}(1)$ exist. From part~2 of 
Proposition~\ref{prop_firstRel2h}
and Corollary~\ref{cor_UVW1vC} we get that
$$
B_h(1)=B_{0,h}(1)+
    C_{0,h}(1)^2B_h(1)\,.
$$
From the proof of part~2 of Theorem~\ref{thm_polya1gen1} we 
know that $B_h(1)\ne0$. Hence $C_{0,h}(1)\in\mathrm{sqrt}\big(1-\frac{B_{0,h}(1)}{B_h(1)}\big)$. 
By part~1 of 
Proposition~\ref{prop_firstRel2h}, Proposition~\ref{prop_onD} and
Proposition~\ref{prop_tauber1vC}, 
$\lim_{n\to\infty}(a_n^h)'=C_{0,h}(1)$.
The stated formula follows. 

3. We assume that $h$ is nonnegative,  
$B_h(1)=+\infty$ and $v\in V(h)$. By Proposition~\ref{prop_onD},  
$D_h(x)=\frac{1}{1-x}$. Hence for every $n\in\N_0$ the coefficients $(a_n^h)'$, 
$b_n^h$, $(b_n^h)'$, $(c_n^h)'$ and $(c_n^h)''$ lie in $[0,1]$ and we can use 
the functions
$$
F_{B_h},\,F_{B_{0,h}},\,F_{C_{0,h}},\,F_{C_{1,h}}\cc[0,\,1)\to[0,\,+\infty)\,.
$$
Since $v\in V(h)$, we have $(c_n^h)''>0$ for some $n>0$ 
(consider the shortest walk in $G(h)$ joining $1$ and $v$). Thus we deduce by 
the same argument as in the proof of the case $d\le2$ of 
Theorem~\ref{thm_explVal2d3} that $F_{B_{0,h}}$ 
is bounded on $[0,1)$.  This 
boundedness, the assumption that 
$B_h(1)=+\infty$, part~2 of 
Proposition~\ref{prop_firstRel2h} and Theorem~\ref{thm_weakAbel} give 
that $C_{0,h}(1)=1$. By part~1 of 
Proposition~\ref{prop_firstRel2h}, $D_h(x)=\frac{1}{1-x}$, and 
Proposition~\ref{prop_tauber1} we get the stated limit that $\lim_{n\to\infty}(a_n^h)'=C_{0,h}(1)=1$.
\eproof

\noindent
Theorems~\ref{thm_polya1gen2n} and \ref{thm_polya1gen2vC} generalize
Theorem~\ref{thm_explVal2d3}.

We extend Corollary~\ref{cor_B0lomenoB} to $\mathrm{W_{\C}}$.

\begin{cor}\label{cor_B0lomenoB2}
Let $v\in\N\setminus\{1\}$ and $h\cc\N_2\to\R_{\ge0}$ be 
a~$v$-transitive light nonnegative weight such that $D_h(x)$ converges on $[0,1)$, the sum $B_h(1)=+\infty$ and $v\in V(h)$. Then
$$
\lim_{x\to1}F_{B_h}(x)=+\infty\,\text{ and }\,\lim_{x\to1}{\textstyle\frac{F_{B_{0,h}}(x)}{F_{B_h}(x)}=0\,.
}
$$
\end{cor}
\proof
The first limit follows from $B_h(1)=+\infty$ and Theorem~\ref{thm_weakAbel}. We know
that on $[0,1)$, $F_{B_h}\ge1$ and
$\frac{F_{B_{0,h}}}{F_{B_h}}=1-F_{C_{0,h}}^2$. Let $x\to1$. Then
$F_{C_{0,h}}\to 1$ by the argument in the proof of the case $d\le2$ of 
Theorem~\ref{thm_explVal2d3} (note that $(c_n^h)''>0$ for some $n>0$ because $v\in V(h)$) and the second limit follows.
\eproof

\section{Concluding remarks}\label{sec_conclR}

Since this paper is already long, we decided to move several concrete
examples in the follow-up paper \cite{klaz_hors2}. There we obtain
analogues of the four main theorems in the previous section for weights
$$
h\cc\N_2\to\C[[x_1,\,x_2,\,\ds,\,x_k]]\,.
$$
The algebras  $\C[[x_1,x_2,\ds,x_k]]$ of formal power series are endowed
with non-Archimedean norms. 
Non-Archimedean infinite series are simpler than the complex ones we 
dealt with above. 
Weighting objects with polynomials and formal power series is a~basic 
technique in enumerative combinatorics.


\begin{thebibliography}{20}

\addcontentsline{toc}{section}{References}

\bibitem{abel}
N.\,H. Abel, Untersuchungen \"uber die Reihe: $1+\frac{m}{1}x+
\frac{m\cdot(m-1)}{1\cdot2}\cdot x^2+\frac{m\cdot(m-1)\cdot(m-2)}{1\cdot2\cdot3}\cdot 
x^3+\ds\,\ds$ u. s. w. {\em Journal f\"ur die reine und angewandte Mathematik} {\bf 1} (1826), 311--339

\bibitem{albe_al}
M.\,H. Albert, Ch. Bean, A. Claesson, \'E. Nadeau, J. Pantone and H. Ulfarsson, Combinatorial Exploration: An algorithmic framework for enumeration, ArXiv:2202.07715v3, 2024, 99 pp.

\bibitem{alex}
G. Alexanderson, {\em The Random Walks of George P\'olya}, The Mathematical Association of America, Washington, DC 2000 

\bibitem{beck}
J. Beck, Recurrence of inhomogeneous random walks, {\em Period. Math. Hung} {\bf 74} (2017), 137--196 

\bibitem{bend_rich}
E.\,A. Bender and L.\,B. Richmond, Correlated random walks, {\em Annals Prob.} {\bf 12} (1984), 274--278

\bibitem{bill}
P. Billingsley, {\em Probability and Measure}. Third Edition, John Wiley $\&$ Sons,  New York 1995

\bibitem{comt}
L. Comtet, {\em Advanced Combinatorics. The Art of Finite and Infinite Expansions}, D.~Reidel, Dordrecht, Holland 1974 

\bibitem{fell}
W. Feller, {\em An Introduction to Probability Theory and Its Applications.} Volume I. Third Edition, John Wiley $\&$ Sons, 
New York 1968

\bibitem{flaj_sedg}
P. Flajolet and R. Sedgewick, {\em Analytic Combinatorics}, Cambridge University Press, Cambridge, UK
2009

\bibitem{fost_good}
F.\,G. Foster and I.\,J. Good, On a~generalization of Polya's random-walk 
theorem, {\em Qart. J. Math. Oxford} {\bf 4} (1953), 120--126 

\bibitem{goul_jack}
I.\,P. Goulden and D.\,M. Jackson, {\em Combinatorial Enumeration}, J. Wiley $\&$ Sons, New York 1983 

\bibitem{grim_wels}
G. Grimmet and D. Welsh, {\em Probability. An Introduction}. Second Edition, Oxford University Press,
Oxford, UK 2014

\bibitem{klaz}
M. Klazar, A~combinatorial model of discrete random walks, in preparation

\bibitem{klaz1}
M. Klazar, Semi-formal symbolic method in enumerative 
combinatorics, in preparation

\bibitem{klaz_hors2}
M. Klazar and R. Horsk\'y, Extending P\'olya's random walker beyond probability~II. Non-Archimedean weights, in preparation


\bibitem{koch}
Y. Kochetkov, An easy proof of Polya’s theorem on random walks, arXiv:1803.00811v1, 2018, 3 pp.

\bibitem{kolm}
A. Kolmogoroff, {\em Grundbegriffe der Wahrscheinlichkeitsrechnung}, Verlag von Julius Springer, Berlin 1933

\bibitem{lang}
K. Lange, Polya’s random walks theorem revisited, {\em Amer. Math. Monthly} {\bf 122} (2015), 1005--1007

\bibitem{lawl}
G.\,F. Lawler, {\em Intersections of Random Walks}, Birkh\"auser, Boston 1991

\bibitem{levi_pere}
D.\,A. Levin and Y. Peres, P\'olya’s theorem on random walks via P\'olya’s urn, {\em Amer. Math. 
Monthly} {\bf 117} (2010), 220--231

\bibitem{nova}
J. Novak, P\'olya’s random walk theorem, {\em Amer. Math. Monthly} {\bf 121} (2014), 711--716

\bibitem{poly}
G. Polya, \"Uber eine Aufgabe der Wahrscheinlichkeitsrechnung betreffend die Irrfahrt
im Strassennetz, {\em Math. Annalen} {\bf 84} (1921), 149--160.


\bibitem{poly_38}
G. Polya, Sur la promenade au hasard dans un r\'eseau de rues, {\em Actualit\'es Sci. Ind.} {\bf 
734} (1938), 25--44

\bibitem{reny}
A. R\'enyi, {\em Teorie pravd\v epodobnosti}, Academia, Praha 1972 [Probability Theory, translation of 
the German edition in 1962, translator not mentioned]

\bibitem{reve}
P. R\'ev\'esz, {\em Random Walk in Random and Nonrandom Environments}, World Scientific Publishing Co., Inc., Teaneck, NJ, 1990

\bibitem{robb}
H. Robbins, A remark on Stirling’s formula, {\em Amer. Math. Monthly} {\bf6} (1955), 26--29





\bibitem{symb}
Symbolic method (combinatorics), Wikipedia article, 
\url{https://en.wikipedia.org/wiki/Symbolic_method_(combinatorics)}

\bibitem{tene}
G. Tenenbaum, {\em Introduction to Analytic and Probabilistic
Number Theory.} Third Edition, AMS, Providence, RI 2015 

\bibitem{wins}
V. Winstein, P\'olya's Theorem on Random Walks, slides, 2021, available at
\url{https://vilas.us/mathnotes/osutalks/ReadingClassics_PolyasTheorem.pdf}

\bibitem{woes}
W. Woess, {\em Random Walks on Infinite Graphs and Groups}, Cambridge University Press, Cambridge, UK
2000

\bibitem{woes_MCh}
W. Woess, {\em Denumerable Markov chains}, 
European Mathematical Society (EMS), Z\"urich 2009

\end{thebibliography}
\end{document}